\documentstyle[11pt]{article}
\begin{document}
%
%
%
\newtheorem{theorem}      {Th\'eor\`eme}[section]
\newtheorem{theorem*}     {theorem}
\newtheorem{proposition}  [theorem]{Proposition}
\newtheorem{definition}   [theorem]{Definition}
\newtheorem{e-lemme}        [theorem]{Lemma}
\newtheorem{cor}   [theorem]{Corollaire}
\newtheorem{resultat}     [theorem]{R\'esultat}
\newtheorem{eexercice}    [theorem]{Exercice}
\newtheorem{rrem}    [theorem]{Remarque}
\newtheorem{pprobleme}    [theorem]{Probl\`eme}
\newtheorem{eexemple}     [theorem]{Exemple}
\newcommand{\preuve}      {\paragraph{Preuve}}
\newenvironment{probleme} {\begin{pprobleme}\rm}{\end{pprobleme}}
\newenvironment{remarque} {\begin{rremarque}\rm}{\end{rremarque}}
\newenvironment{exercice} {\begin{eexercice}\rm}{\end{eexercice}}
\newenvironment{exemple}  {\begin{eexemple}\rm}{\end{eexemple}}
%
%
\newtheorem{e-theo}      [theorem]{Theorem}
\newtheorem{theo*}     [theorem]{Theorem}
\newtheorem{e-pro}  [theorem]{Proposition}
\newtheorem{e-def}   [theorem]{Definition}
\newtheorem{e-lem}        [theorem]{Lemma}
\newtheorem{e-cor}   [theorem]{Corollary}
\newtheorem{e-resultat}     [theorem]{Result}
\newtheorem{ex}    [theorem]{Exercise}
\newtheorem{e-rem}    [theorem]{Remark}
\newtheorem{prob}    [theorem]{Problem}
\newtheorem{example}     [theorem]{Example}
\newcommand{\proof}         {\paragraph{Proof~: }}
\newcommand{\hint}          {\paragraph{Hint}}
\newcommand{\heuristicproof}{\paragraph{heuristic proof}}
\newenvironment{e-probleme} {\begin{e-pprobleme}\rm}{\end{e-pprobleme}}
\newenvironment{e-remarque} {\begin{e-rremarque}\rm}{\end{e-rremarque}}
\newenvironment{e-exercice} {\begin{e-eexercice}\rm}{\end{e-eexercice}}
\newenvironment{e-exemple}  {\begin{e-eexemple}\rm}{\end{e-eexemple}}
\newcommand{\reell}    {{{\rm I\! R}^l}}
\newcommand{\reeln}    {{{\rm I\! R}^n}}
\newcommand{\reelk}    {{{\rm I\! R}^k}}
\newcommand{\reelm}    {{{\rm I\! R}^m}}
\newcommand{\reelp}    {{{\rm I\! R}^p}}
\newcommand{\reeld}    {{{\rm I\! R}^d}}
\newcommand{\reeldd}   {{{\rm I\! R}^{d\times d}}}
\newcommand{\reelnn}   {{{\rm I\! R}^{n\times n}}}
\newcommand{\reelnd}   {{{\rm I\! R}^{n\times d}}}
\newcommand{\reeldn}   {{{\rm I\! R}^{d\times n}}}
\newcommand{\reelkd}   {{{\rm I\! R}^{k\times d}}}
\newcommand{\reelkl}   {{{\rm I\! R}^{k\times l}}}
\newcommand{\reelN}    {{{\rm I\! R}^N}}
\newcommand{\reelM}    {{{\rm I\! R}^M}}
\newcommand{\reelplus} {{{\rm I\! R}^+}}
\newcommand{\reelo}    {{{\rm I\! R}\setminus\{0\}}}
\newcommand{\reld}    {{{\rm I\! R}_d}}
\newcommand{\relplus} {{{\rm I\! R}_+}}
\newcommand{\1}        {{\bf 1}}

\newcommand{\cov}      {{\hbox{cov}}}
\newcommand{\sss}      {{\cal S}}
\newcommand{\indic}    {{{\rm I\!\! I}}}
\newcommand{\pp}       {{{\rm I\!\!\! P}}}
\newcommand{\qq}       {{{\rm I\!\!\! Q}}}
\newcommand{\ee}       {{{\rm I\! E}}}

\newcommand{\B}        {{{\rm I\! B}}}
\newcommand{\cc}       {{{\rm I\!\!\! C}}}
\newcommand{\HHH}        {{{\rm I\! H}}}
\newcommand{\N}        {{{\rm I\! N}}}
\newcommand{\R}        {{{\rm I\! R}}}
\newcommand{\D}        {{{\rm I\! D}}}
\newcommand{\Z}       {{{\rm Z\!\! Z}}}
\newcommand{\C}        {{\bf C}}        
\newcommand{\T}        {{\bf T}}        
\newcommand{\E}        {{\bf E}}        
\newcommand{\rfr}[1]    {\stackrel{\circ}{#1}}
\newcommand{\equiva}    {\displaystyle\mathop{\simeq}}
\newcommand{\eqdef}     {\stackrel{\triangle}{=}}
\newcommand{\limps}     {\mathop{\hbox{\rm lim--p.s.}}}
\newcommand{\Limsup}    {\mathop{\overline{\rm lim}}}
\newcommand{\Liminf}    {\mathop{\underline{\rm lim}}}
\newcommand{\Inf}       {\mathop{\rm Inf}}
\newcommand{\vers}      {\mathop{\;{\rightarrow}\;}}
\newcommand{\versup}    {\mathop{\;{\nearrow}\;}}
\newcommand{\versdown}  {\mathop{\;{\searrow}\;}}
\newcommand{\vvers}     {\mathop{\;{\longrightarrow}\;}}
\newcommand{\cvetroite} {\mathop{\;{\Longrightarrow}\;}}
\newcommand{\ieme}      {\hbox{i}^{\hbox{\smalltype\`eme}}}
\newcommand{\eqps}      {\, \buildrel \rm \hbox{\rm\smalltype p.s.} \over =
\,}
\newcommand{\eqas}      {\,\buildrel\rm\hbox{\rm\smalltype a.s.} \over = \,}
\newcommand{\argmax}    {\hbox{{\rm Arg}}\max}
\newcommand{\argmin}    {\hbox{{\rm Arg}}\min}
\newcommand{\indep}{\perp\!\!\!\!\perp}
\newcommand{\abs}[1]{\left| #1 \right|}
\newcommand{\crochet}[2]{\langle #1 \,,\, #2 \rangle}
\newcommand{\espc}[3]   {E_{#1}\left(\left. #2 \right| #3 \right)}
\newcommand{\rang}{\hbox{rang}}
\newcommand{\rank}{\hbox{rank}}
\newcommand{\signe}{\hbox{signe}}
\newcommand{\sign}{\hbox{sign}}

\newcommand\hA{{\widehat A}}
\newcommand\hB{{\widehat B}}
\newcommand\hC{{\widehat C}}
\newcommand\hD{{\widehat D}}
\newcommand\hE{{\widehat E}}
\newcommand\hF{{\widehat F}}
\newcommand\hG{{\widehat G}}
\newcommand\hH{{\widehat H}}
\newcommand\hI{{\widehat I}}
\newcommand\hJ{{\widehat J}}
\newcommand\hK{{\widehat K}}
\newcommand\hL{{\widehat L}}
\newcommand\hM{{\widehat M}}
\newcommand\hN{{\widehat N}}
\newcommand\hO{{\widehat O}}
\newcommand\hP{{\widehat P}}
\newcommand\hQ{{\widehat Q}}
\newcommand\hR{{\widehat R}}
\newcommand\hS{{\widehat S}}
\newcommand\hTT{{\widehat T}}
\newcommand\hU{{\widehat U}}
\newcommand\hV{{\widehat V}}
\newcommand\hW{{\widehat W}}
\newcommand\hX{{\widehat X}}
\newcommand\hY{{\widehat Y}}
\newcommand\hZ{{\widehat Z}}

\newcommand\ha{{\widehat a}}
\newcommand\hb{{\widehat b}}
\newcommand\hc{{\widehat c}}
\newcommand\hd{{\widehat d}}
\newcommand\he{{\widehat e}}
\newcommand\hf{{\widehat f}}
\newcommand\hg{{\widehat g}}
\newcommand\hh{{\widehat h}}
\newcommand\hi{{\widehat i}}
\newcommand\hj{{\widehat j}}
\newcommand\hk{{\widehat k}}
\newcommand\hl{{\widehat l}}
\newcommand\hm{{\widehat m}}
\newcommand\hn{{\widehat n}}
\newcommand\ho{{\widehat o}}
\newcommand\hp{{\widehat p}}
\newcommand\hq{{\widehat q}}
\newcommand\hr{{\widehat r}}
\newcommand\hs{{\widehat s}}
\newcommand\htt{{\widehat t}}
\newcommand\hu{{\widehat u}}
\newcommand\hv{{\widehat v}}
\newcommand\hw{{\widehat w}}
\newcommand\hx{{\widehat x}}
\newcommand\hy{{\widehat y}}
\newcommand\hz{{\widehat z}}

\newcommand\tA{{\widetilde A}}
\newcommand\tB{{\widetilde B}}
\newcommand\tC{{\widetilde C}}
\newcommand\tD{{\widetilde D}}
\newcommand\tE{{\widetilde E}}
\newcommand\tF{{\widetilde F}}
\newcommand\tG{{\widetilde G}}
\newcommand\tH{{\widetilde H}}
\newcommand\tI{{\widetilde I}}
\newcommand\tJ{{\widetilde J}}
\newcommand\tK{{\widetilde K}}
\newcommand\tL{{\widetilde L}}
\newcommand\tM{{\widetilde M}}
\newcommand\tN{{\widetilde N}}
\newcommand\tOO{{\widetilde O}}
\newcommand\tP{{\widetilde P}}
\newcommand\tQ{{\widetilde Q}}
\newcommand\tR{{\widetilde R}}
\newcommand\tS{{\widetilde S}}
\newcommand\tTT{{\widetilde T}}
\newcommand\tU{{\widetilde U}}
\newcommand\tV{{\widetilde V}}
\newcommand\tW{{\widetilde W}}
\newcommand\tX{{\widetilde X}}
\newcommand\tY{{\widetilde Y}}
\newcommand\tZ{{\widetilde Z}}

\newcommand\ta{{\widetilde a}}
\newcommand\tb{{\widetilde b}}
\newcommand\tc{{\widetilde c}}
\newcommand\td{{\widetilde d}}
\newcommand\te{{\widetilde e}}
\newcommand\tf{{\widetilde f}}
\newcommand\tg{{\widetilde g}}
\newcommand\th{{\widetilde h}}
\newcommand\ti{{\widetilde i}}
\newcommand\tj{{\widetilde j}}
\newcommand\tk{{\widetilde k}}
\newcommand\tl{{\widetilde l}}
\newcommand\tm{{\widetilde m}}
\newcommand\tn{{\widetilde n}}
\newcommand\tio{{\widetilde o}}
\newcommand\tp{{\widetilde p}}
\newcommand\tq{{\widetilde q}}
\newcommand\tr{{\widetilde r}}
\newcommand\ts{{\widetilde s}}
\newcommand\tit{{\widetilde t}}
\newcommand\tu{{\widetilde u}}
\newcommand\tv{{\widetilde v}}
\newcommand\tw{{\widetilde w}}
\newcommand\tx{{\widetilde x}}
\newcommand\ty{{\widetilde y}}
\newcommand\tz{{\widetilde z}}

\newcommand\bA{{\overline A}}
\newcommand\bB{{\overline B}}
\newcommand\bC{{\overline C}}
\newcommand\bD{{\overline D}}
\newcommand\bE{{\overline E}}
\newcommand\bFF{{\overline F}}
\newcommand\bG{{\overline G}}
\newcommand\bH{{\overline H}}
\newcommand\bI{{\overline I}}
\newcommand\bJ{{\overline J}}
\newcommand\bK{{\overline K}}
\newcommand\bL{{\overline L}}
\newcommand\bM{{\overline M}}
\newcommand\bN{{\overline N}}
\newcommand\bO{{\overline O}}
\newcommand\bP{{\overline P}}
\newcommand\bQ{{\overline Q}}
\newcommand\bR{{\overline R}}
\newcommand\bS{{\overline S}}
\newcommand\bT{{\overline T}}
\newcommand\bU{{\overline U}}
\newcommand\bV{{\overline V}}
\newcommand\bW{{\overline W}}
\newcommand\bX{{\overline X}}
\newcommand\bY{{\overline Y}}
\newcommand\bZ{{\overline Z}}

\newcommand\ba{{\overline a}}
\newcommand\bb{{\overline b}}
\newcommand\bc{{\overline c}}
\newcommand\bd{{\overline d}}
\newcommand\be{{\overline e}}
\newcommand\bff{{\overline f}}
\newcommand\bg{{\overline g}}
\newcommand\bh{{\overline h}}
\newcommand\bi{{\overline i}}
\newcommand\bj{{\overline j}}
\newcommand\bk{{\overline k}}
\newcommand\bl{{\overline l}}
\newcommand\bm{{\overline m}}
\newcommand\bn{{\overline n}}
\newcommand\bo{{\overline o}}
\newcommand\bp{{\overline p}}
\newcommand\bq{{\overline q}}
\newcommand\br{{\overline r}}
\newcommand\bs{{\overline s}}
\newcommand\bt{{\overline t}}
\newcommand\bu{{\overline u}}
\newcommand\bv{{\overline v}}
\newcommand\bw{{\overline w}}
\newcommand\bx{{\overline x}}
\newcommand\by{{\overline y}}
\newcommand\bz{{\overline z}}

%
\newcommand{\AAA}{{\cal A}}
\newcommand{\BB}{{\cal B}}
\newcommand{\CC}{{\cal C}}
\newcommand{\DD}{{\cal D}}
\newcommand{\EE}{{\cal E}}
\newcommand{\FF}{{\cal F}}
\newcommand{\GG}{{\cal G}}
\newcommand{\HH}{{\cal H}}
\newcommand{\II}{{\cal I}}
\newcommand{\JJ}{{\cal J}}
\newcommand{\KK}{{\cal K}}
\newcommand{\LL}{{\cal L}}
\newcommand{\NN}{{\cal N}}
\newcommand{\MM}{{\cal M}}
\newcommand{\OO}{{\cal O}}
\newcommand{\PP}{{\cal P}}
\newcommand{\QQ}{{\cal Q}}
\newcommand{\RR}{{\cal R}}
\newcommand{\SS}{{\cal S}}
\newcommand{\TT}{{\cal T}}
\newcommand{\UU}{{\cal U}}
\newcommand{\VV}{{\cal V}}
\newcommand{\WW}{{\cal W}}
\newcommand{\XX}{{\cal X}}
\newcommand{\YY}{{\cal Y}}
\newcommand{\ZZ}{{\cal Z}}
\newcommand{\tbullet}{$\bullet$}
\newcommand{\ot}{\leftarrow}
\newcommand{\carre}{\hfill$\Box$}
\newcommand{\carreb}{\hfill\rule{0.25cm}{0.25cm}}
%
%
\newcommand{\dontforget}[1]
{{\mbox{}\\\noindent\rule{1cm}{2mm}\hfill don't forget : #1
\hfill\rule{1cm}{2mm}}\typeout{---------- don't forget : #1 ------------}}
\newcommand{\note}[2]
{ \noindent{\sf #1 \hfill \today}

\noindent\mbox{}\hrulefill\mbox{}
\begin{quote}\begin{quote}\sf #2\end{quote}\end{quote}
\noindent\mbox{}\hrulefill\mbox{}
\vspace{1cm}
}
\newcommand{\rond}[1]     {\stackrel{\circ}{#1}}
\newcommand{\rondf}       {\stackrel{\circ}{\FF}}
\newcommand{\point}[1]     {\stackrel{\cdot}{#1}}

\newcommand\relatif{{\rm \rlap Z\kern 3pt Z}}

\title{\huge   Segre varieties, CR geometry and Lie
symmetries of second order PDE systems.  }
\author{Alexandre Sukhov }
\date{}
\maketitle

Abstract. We establish a link between the CR geometry of real 
analytic submanifolds in $\cc^n$ and the geometric PDE theory.
The main idea of our approach is to consider biholomorphisms of 
a Levi-nondegenerate real analytic Cauchy-Riemann manifold ${\cal M}$ as
poinwise symmetries of a second order holomorphic PDE system 
defining the Segre family of ${\cal M}$. This allows to employ the 
well-elaborated PDE tools in order to study the biholomorphism 
group of ${\cal M}$. We give several examples and applications to the 
CR geometry: the results on the finite dimensionality of the 
biholomorphism group and precise estimates of its dimension, 
explicit parametrization of the Lie algebra of infinitesimal 
automorphisms etc. We deduce these results as a special case 
of more general statements concerning related properties of 
symmetries of second order PDE systems.



\section{ Introduction}

The main goal of the present paper is to establish the relationship between
the
{\cal CR} geometry of a real analytic generic submanifold of $\cc^{n+m}$
and the geometric (or formal) theory of PDE. We apply a general method
which is due to S.Lie in order to study infinitesimal symmetries of a
holomorphic completely overdetermined involutive second order PDE system
with first order relations and $n$ independent and $m$ dependent variables.
For any given system of this class this method allows to determine whether
the
dimension of the Lie algebra of infinitesimal symmetries if finite; if this
is the
case, the Lie method leads to
  explicit recurcive formulae which permit to
compute terms of any order in the Taylor expansion of coefficients of an
infinitesimal symmetry of such a system and to  show that these expansions
(and so any symmetry) are uniquely determined by their  terms of  finite
order.
This gives a
 precise upper estimate of the dimension of the symmetry group for such a
system
and an explicit parametrization of the symmetry group.

From the complex analysis point of view our interest in these questions 
is explained by the fact that  the Segre family of a real analytic generic Levi
nondegenerate subvariety ${\cal M}$ in $\cc^{n+m}$ (introduced to 
the modern theory by S.S.Chern and S.Webster) is
a family of (graphs of) solutions of a holomorphic completely overdetermined
involutive PDE system with  $m$  dependent and $n$ independent variables and
some additional first order equations if the real codimension of ${\cal M}$
is
$> 1$ , i.e. if ${\cal M}$ is not a hypersurface. Systems without first
order
relations were studied in our previous paper  \cite{Su} so in the present
paper we consider the  more complicated higher codimensional
case.
The biholomorphic invariance of the Segre family means precisely that every
biholomorphic automorphism of ${\cal M}$ is a Lie symmetry of the PDE system
defining its Segre family i.e. maps the graph of a
solution to the graph of another solution. So we show how PDE symmetries techniques
can be used in order to  study the    complex geometry of real
analytic submanifolds in $\cc^n$ and to obtain
  precise upper estimates of the dimension
 and explicit parametrization  of their automorphism groups
etc.; various results of this type have been obtained by several authors 
using different methods
(see a more detailed discussion below).
 But it is worth to emphasize that systems describing the Segre
 families of real analytic submanifolds form {\it a very special subclass}
 in the class of holomorphic completely overdetermined involutive systems
with first order relations.
So we consider a much more general situation and generalize some
 known results on automorphisms of CR manifolds.

In the present paper we pay more attention to the development of basic 
tools of the proposed PDE approach to the CR geometry and do not consider 
the most general classes of CR manifolds in order to avoid technical 
complications and long computations. However, the proposed method allows
to obtain much more general and precise results not only for CR 
manifolds, but for symmetries of wide classes of PDE as well. 
 Our main conclusion is that the very intesively  developping theory of
CR maps can be naturally viewed as a part of the  geometric PDE theory and
actually studies special poinwise symmetries of special holomorphic PDE
systems. From
our point of view, the further progress in the study of CR maps between real
analytic submanifolds in $\cc^n$ may be achieved by application of advanced
tools of the formal PDE theory combining with complex algebraic and
differential geometry methods.
This provides the natural framework for the CR geometry of real analytic
manifolds  and links it with the classical complex geometry.

\section{ Generatities of the Lie theory}

In this section we recall certain basic tools of the Lie method of study
of infinitesimal symmetries of differential equations. They are very well
known to
the experts in the geometric PDE theory and the differential geometry; for
reader's convenient we give a brief exposition. A more detailed information
and the proofs
of all statements of this section can be found in \cite{BlKu}, \cite{Ol},
\cite{Ov},
\cite{Po}.

{\bf 2.1. Local transformation groups and symmetry groups.}

 Let $\Omega$ be a domain in
$\cc^n$. A local group of biholomorphic transformations
acting on $\Omega$ is given by a (local) connected complex Lie group $G$, a
domain
$D$ such that $\{ e \} \times \Omega \subset D \subset G \times \Omega$, and
a
holomorphic map $\varphi: D \longrightarrow \Omega$ with the following
properties:
(i) if $(h,x) \in D$, $(g,\varphi(h,x)) \in D$, and also $(gh,x) \in D$,
then
$\varphi(g,\varphi(h,x)) = \varphi(gh,x)$; (ii) for all $x \in \Omega$,
$\varphi(e,x) =x$; (iii) if $(g,x) \in D$, then $(g^{-1}, \varphi(g,x)) \in
D$ and
$\varphi(g^{-1},\varphi(g,x)) = x$.

Historically the notion of a group of transformations was introduced by
S.Lie
in connection with a study of transformations preserving a given PDE system
(or more precisely, the space of its solutions). Such transformations are
called
{\it symmetries} (sometimes, the Lie symmetries, pointwise
symmetries, classical symmetries). In the present paper we apply the Lie
method of studying of
PDE symmetries to a special but geometrically important class of holomorphic
completely overdertermined
second order PDE systems with first order relations, i.e. systems of the form

\begin{eqnarray*}
& &({\cal S}): u^1_{x_ix_j} = F_{ij}(x,u,u^1_x),  i,j =
1,...,n,\\
& &u^k_x = H^k(x,u,u^1_x), k = 2,...,m
\end{eqnarray*}
where $x = (x_1,...,x_n)$ are independent variables, $u(x) =
(u^1(x),...,u^m(x))$
are unknown functions (dependent variables), $u^j_x = (u^j_{x_1},...,
u^j_{x_n})$ and
$F_{ij}$, 
$H^k$ are holomorphic functions (of course, we will always assume that $F_{ij} = 
F_{ji}$). Since this system is highly 
overdetermined, it is natural to assume that it satisfies some 
compatibility conditions. We will  assume that such a system
satisfies some
integrability conditions of the Frobenius type (see below).
 This class of systems naturally arises
in various areas
of the geometry and PDE.

The solutions of such a system are holomorphic vector valued functions $u =
u(x)$;
denote by $\Gamma_u$ the graph of a solution $u$.

\begin{e-def}
A symmetry group $Sym(S)$ of a system $({\cal S})$ is a local complex
transformation
group acting on a domain in the space $\cc^n_x \times \cc^m_u$ of independent and
dependent variables
with the following property: for every solution $u(x)$ of $({\cal S})$ and
every $g \in G$ such
that the image $g(\Gamma_u)$ is defined, it is a graph of a solution of
$({\cal S})$.
\end{e-def}

Often the {\it largest} symmetry group is of main interest (and so
we write {\it the} symmetry group); for us this is not very
essential since our methods give a description of ${\it any}$ symmetry group
for given system. In order to fix the terminology, {\it everywhere below 
by the symmetry group we mean the largest one}.

The definition of a symmetry group given above is not very well working in
practice
in the sense that it does not give an efficient tool to find the Lie
symmetries. The
main idea of the Lie method is to study the Lie algebra of a symmetry group
instead of the group
itself.

{\bf 2.2. Vector fields and one-parameter transformation groups.} Consider a
one parameter
local complex Lie group of transformations (LTG) $x^{*} = X(x,t)$ with the
identity
$t = 0$ acting on a complex manifold with local coordinates $x =
(x_1,...,x_n)$.
Let $\theta(x) = \frac{\partial X}{\partial t}\vert_{t = 0}$.
The vector field $X = X(x) = \sum_{j=1}^n
\theta_j(x)\frac{\partial}{\partial x_j}$
is called the {\it infinitesimal generator} of our (LTG) ; we use the vector
notation: if $f = (f_1,...,f_k)$ is a holomorphic vector function,
then $Xf = (Xf_1,...,
Xf_k)$. In particular, $Xx = \theta(x)$. Recall that there exists a
parametrization
$\tau(t)$ such that the above (LTG)  is equivalent to the solution of the
initial
value problem for the first order ODE system $\frac{dx^*}{d\tau} =
\theta(x^*)$
(the First Fundamental Lie Theorem). A one-parameter (LTG) can be found from
its
infinitesimal
generator by means of the Lie series (the exponential map):
$x^* = e^{tX}x = x + tXx + (t^2/2)X^2x + ...$, where $X^k:= X X^{k-1}$,
$k = 1,2,...$, $X^{0}f(x):= f(x)$, $t \in \cc$. In the general case of a
$d$-
dimensional Lie transformation group $G$ any group element in a neighborhood
of the
identity can be obtained by the exponential map for a suitable vector field
from the
Lie algebra of $G$. So every local Lie group is completely determined by a
vector
field basis
$\{ X_1,...,X_d\}$ of its Lie algebra and can be explicitely parametrized
via the
exponential map $e^{\sum t_jX_j} = \Pi e^{t_jX_j}$; the parameters
$t_1,...,t_d$ are
local coordinates on $G$. The exponential map can be used as a definition of
a
symmetry group; this group is a finite dimensional Lie group if and only its
Lie algebra is finite dimensional.

{\bf 2.3. Jet bundles and prolongations of group actions.} The second key
tool of
the Lie theory is the notion of {\it prolongation} of an LTG action to a jet
bundle.
Recall this construction. Let $f$ and $g$ be two holomorphic maps in a
neighborhood of the
origin in $\cc^n$ to $\cc^m$ taking the origin to the origin. As usual, we
say that they
have the same $r$-jet at the origin if $(\partial^{\alpha}f)(0) =
(\partial^{\alpha}g)(0)$
for every $\alpha: \vert \alpha \vert \leq r$ where we use the
following notation (which we will keep everywhere through this paper): 
$\partial^{\alpha} \varphi = \frac{\partial^r \varphi}{\partial
x_{\alpha_1}...
\partial x_{\alpha_r}}$
for $\alpha = (\alpha_1,...,\alpha_r)$, $\alpha_1 \leq ...\leq
\alpha_r$, and $\vert \alpha \vert := r$.

 More generally,  let $M$ and
$N$ be two
complex manifolds and $f: M \longrightarrow N$, $g:M \longrightarrow N$ be
two
holomorphic maps. Let $x$ and $u$ be local holomorphic coordinates near $p
\in M$ and
$q \in N$ respectively such that $x(p) = 0$, $u(q) = 0$. We say that $f$ and
$g$ have the
same $r$-jet at $p$, if $u \circ f \circ x^{-1}$ and $u \circ g \circ
x^{-1}$ have the
same $r$-jet.  It is easy to see that the definition is correct, i.e. does
not depend
on the choice of the coordinates. The relation that two maps have the same
$r$-jet
at $p$ is an equivalence relation and the equivalence class with the
representative $f$
is denoted by $j^r_p(f)$; it is called the $r$- {\it jet} of $f$ at $p$.
The point $p$ is called the source and the point $q$ the target
of $j^q_p(f)$. Denote by $J_{p,q}^r(M,N)$ the set of all $r$-jets of maps
from $M$ to $N$ with the source $p$ and the target $q$ and consider the set
$J^r(M,N) = U_{p \in M, q \in N} J^r_{p,q}(M,N)$.
Consider also the natural projections $\pi_M: J^r(M,N) \longrightarrow M$
and
$\pi_N: J^r(M,N) \longrightarrow N$ defined by $\pi_M(j_p^r(f)) = p$ and
$\pi_N(j_p^r(f)) = f(p)$. Declaring the pullbacks of open sets in $M$ and
$N$
to be open, we define the natural topology on $J^r(M,N)$. Using local
coordinates
$x$ on $M$ and $u$ on $N$ defined as above we may define a local coordinate
system $h$ on
$J^r(M,N)$ as follows. Set $u^{(1)} =
(u^1_1,...,u^1_n,...,u^{m}_1,...,u^{m}_n)$,...
, $u^{(s)} = (u^j_{\alpha})$ with $j= 1,...,m$, $\alpha =
(\alpha_1,...\alpha_s)$,
$\alpha_1\leq \alpha_2 \leq ...\leq \alpha_s$. The chart $h$ is defined by

\begin{eqnarray*}
& &h: j^r_{p'}(f) \mapsto (x_j,u^k, u^{(1)},...,u^{(r)})\\
& &x_j = x_j(p'), u^k =u^{k}(f(p')), u^{j}_{\alpha_1...\alpha_s} =
\partial^{\alpha}
(u^j \circ f \circ x^{-1})(x(p')), 1 \leq s \leq r
\end{eqnarray*}
for $p'$ close enough to $p$.
 These coordinates are called the {\it natural coordinates} on $J^r(M, N)$.
The Leibnitz formula and
the chain rule imply that  biholomorphic changes of local coordinates on $M$
and $N$ induce a biholomorphic change of local
coordinates in $J^r(M, N)$. This defines the natural structure of a complex
manifold on the space $J^r(M,N)$ and equips it with the structure of a
holomorphic fiber bundle over $M \times N$ with the natural projection 
$\pi_{M \times N}: J^r(M,N) \longrightarrow M \times N$.

Let $G$ be a local group of biholomorphic transformations acting on $M \times
N$.
Every  biholomorphism $g \in G$ , $g: M \times N \longrightarrow M \times
N$,
$g: (x,u) \mapsto (x^*,u^*)$ {\it close enough to the identity}
lifts canonically to a fiber preserving biholomorphism $g^{(r)}: J^r(M,N)
\longrightarrow J^r(M,N)$ as follows: if $u = f(x)$ is a holomorphic
function
near $p$, $q = f(p)$ and $u^* = f^*(x^*)$ is its image under $g$ (that is the
graph of $f^*$ is the
image of the graph of $f$ under $g$ near the point $(p^*, q^*) = g(p,q)$),
then the jet $j^r_{p^*}(f^*)$ is by the definition the image of $j^r_p(f)$
under
$g^{(r)}$. In particular, a one-parameter local Lie group of
transformations $G$ canonically lifts to $J^r(M,N)$ as a one-parameter Lie
group
of transformations $G^{(r)}$  which is called the {\it r-prolongation} of
$G$. The infinitesimal
generator $X^{(r)}$ of $G$ is called the {\it r-prolongation} of the
infinitesimal
generator $X$ of $G$.

Our considerations will be purely local so $M$ and $N$ will be open 
subsets in $\cc^n$ and $\cc^m$ respectively. In this case we write 
$J^1(n,m)$ instead of $J^1(M,N)$.

Consider in local coordinates a vector field $X(x,u) = \sum_{j=1}^{n}
\theta^j(x,u)\frac{\partial}{\partial x_j} + \sum_{k=1}^m
\eta^{k}(x,u)\frac{\partial}{\partial u^k}$. In the natural coordinates
its $r$-prolongation  has the form

\begin{eqnarray*}
X^{(r)} = X + \sum_{j,\mu} \eta^{\mu}_j\frac{\partial}{\partial u^{\mu}_j}
+...
+ \sum_{i_1,...,i_r,\mu}
\eta^{\mu}_{i_1i_2...i_r}\frac
{\partial}{\partial u^{\mu}_{i_1i_2...i_r}}
\end{eqnarray*}

In order to compute the coefficients of this prolongation, define the
operator of
{\it total derivative}:

\begin{eqnarray*}
D_i = \frac{\partial}{\partial x_i} + \sum_k u^k_i \frac{\partial}{\partial
u^k} +
\sum_{\mu,j} u^{\mu}_{ij}\frac{\partial}{\partial u^{\mu}_j} + ...
\end{eqnarray*}

The following elementary statement gives an explicit recursive formula for
the coefficients
of a prolongation and is the main computational tool in the Lie theory.

\begin{e-pro}
\label{pro2.3}
One has

\begin{eqnarray*}
\eta^{\mu}_i = D_i\eta^{\mu} - \sum_j (D_i\theta^j) u^{\mu}_j,
\eta^{\mu}_{i_1...i_{r-1}i_r} = D_{i_r} \eta^{\mu}_{i_1...i_{r-1}} -
\sum_j (D_{i_r}\theta^j)u^{\mu}_{i_1...i_{r-1}j}
\end{eqnarray*}
\end{e-pro}

In particular  the second prolongation $X^{(2)}$ is given by $X^{(2)} = X^{(1)} + \sum_{\mu; i_1 \neq i_2}
\eta^{\mu}_{i_1i_2}\frac{\partial}
{\partial u^{\mu}_{i_1i_2}} + \sum_{\mu;i} 
\eta^{\mu}_{ii}\frac{\partial}{\partial u^{\mu}_{ii}}$ with $X^{(1)} =X + \sum_{\mu,i} \eta^{\mu}_i\frac{\partial}{\partial
u^{\mu}_i}$

 Proposition \ref{pro2.3} implies the following formula giving an 
explicit
expression for the coefficients of $X^{(2)}$:

\begin{eqnarray*}
& &\eta^{\mu}_{i_1} = \eta^{\mu}_{x_{i_1}} + 
\sum_k u^k_{i_1}
 \eta^{\mu}_{ u^k} - \sum_j \left (  
\theta^j_{ x_{i_1}} + \sum_k u^k_{i_1} \theta^j_{ 
u^k} \right )
u^{\mu}_j ,\\
& &\eta^{\mu}_{i_1i_2} = \eta^{\mu}_{ x_{i_2}
 x_{i_1}} +u^{\mu}_{i_1} \left [\eta^{\mu}_{ x_{i_2} u^{\mu}} 
- \theta^{i_1}_{ x_{i_2} x_{i_1}} \right ] + u^{\mu}_{i_2}
\left [\eta^{\mu}_{ x_{i_1} u^{\mu}}- 
\theta^{i_2}_{ x_{i_2} x_{i_1}}\right ]+ 
\sum_{k \neq \mu} u^k_{i_1}\eta^{\mu}_{ x_{i_2} u^k}\\
& & + \sum_{k\neq \mu} u^k_{i_2}\eta^{\mu}_
{x_{i_1}u^k} -
\sum_{k\neq i_1,k\neq i_2}u^{\mu}_k\theta^k_{ x_{i_2} x_{i_1}}-
\sum_{k; j\neq i_2}u^k_{i_1}u^{\mu}_{j}\theta^j_{x_{i_2} u^k}-
\sum_{i; s\neq i_1} u^i_{i_2}u^{\mu}_{s}\theta^s_{ x_{i_1} u^i}\\
& & + \sum_{r\neq \mu, p\neq \mu}u^r_{i_2}u^{p}_{i_1}
 \eta^{\mu}_{ u^r  u^p} +
\sum_{t \neq \mu} u^t_{i_1}u^{\mu}_{i_2}\left 
[-\theta^{i_2}_{ x_{i_2} u^t} +
\eta^{\mu}_{ u^{\mu}  u^t}\right ] +
\sum_{q \neq \mu} u^q_{i_2}u^{\mu}_{i_1}\left 
[-\theta^{i_1}_{ u^q  x_{i_1}} + \eta^{\mu}_{ u^q  u^{\mu}}\right ] \\
& &+ \left [\eta^{\mu}_{ u^{\mu} u^{\mu}} - \theta^{i_2}_{ x_{i_2} u^{\mu}}
 -  \theta^{i_1}_{ x_{i_1}  u^{\mu}} \right ]
u^{\mu}_{i_1}u^{\mu}_{i_2} -
\sum_{a,b,s} u^a_{i_2}u^b_{i_1}u^{\mu}_{s}\theta^s_
{ u^a  u^b} + \Lambda^{\mu}_{i_1i_2}
\end{eqnarray*}
for $i_1 \neq i_2$ and

\begin{eqnarray*}
& & \eta^{\mu}_{ii} = \eta^{\mu}_{ x_i x_i} +
u^{\mu}_i\left [2\eta^{\mu}_{ x_i u^{\mu}} - 
\theta^{i}_{x_i x_i}\right ] +
2\sum_{k \neq \mu}u^k_i\eta^{\mu}_{x_i
u^k}  -
\sum_{k \neq i}u^{\mu}_k\theta^k_{ x_i x_i} -
2\sum_{k; j \neq i}
u^k_i u^{\mu}_j\theta^j_{ x_i  u^k}\\
& & +
\sum_{r\neq \mu; p\neq \mu} 
u^r_i u^p_i\eta^{\mu}_{ u^r u^p} +
\sum_{t \neq \mu}u^t_iu^{\mu}_i\left [- \theta^i_{ x_i  u^t} +
\eta^{\mu}_{ u^{\mu}  u^t}\right ] +
\sum_{q \neq \mu} u^q_iu^{\mu}_i
\left [-\theta^i_{ x_i  u^q} +
\eta^{\mu}_{ u^q  u^{\mu}}\right ]\\
& &+ \left [\eta^{\mu}_{ u^{\mu} u^{\mu}} -
2 \theta^i_
{ x_i  u^{\mu}}\right ](u^{\mu}_i)^2 -
 \sum_{a,b,s}u^a_iu^b_iu^{\mu}_s \theta^s_{ u^a  u^b} + \Lambda^{\mu}_{ii}
\end{eqnarray*}

with

\begin{eqnarray*}
& &\Lambda^{\mu}_{i_1i_2} = \sum_s u^s_{i_2i_1} \eta^{\mu}_{ u^s} -
\sum_p u^{\mu}_{i_2p}  \theta^p_{ x_{i_1}} - 
\sum_j u^{\mu}_{i_1j}
 \theta^j_{ x_{i_2}} - 
\sum_{p,q} u^q_{i_2i_1} u^{\mu}_p  \theta^p_{ u^q} -
\sum_{p,q} u^{\mu}_{i_2p}u^q_{i_1}  \theta^p_{ u^q} -
\sum_{s,j} u^{\mu}_{i_1j}u^s_{i_2} \theta^j_{ u^s}
 \end{eqnarray*}

{\bf 2.4. Infinitesimal symmetries of differential equations.} An
infinitesimal
generator of a one-parameter group of symmetries of a system of PDE $({\cal
S})$ is called
{\it an infinitesimal symmetry} of this system. They form a  Lie algebra
with
respect to the Lie bracket which  is  denoted by $Lie({\cal S})$.

Let $({\cal S})$ be
a holomorphic PDE system of  order $r$; we consider its
solutions on $M$ with values in $N$. Then it defines naturally a complex
subvariety
$({\cal S}_{r})$ in the jet space $J^{(r)}(M,N)$ obtained by the replacing
of
the derivatives of  dependent variables  by the corresponding 
natural coordinates in the jet space.

{\bf Example 1.} Let $M \subset \cc^2$, $N \subset \cc$ be domains, $({\cal
S})$
be a holomorphic second order ODE $u_{xx} = F(x,u,u_x)$, $(x,u) \in M \times
N$. Let $(x,u,u_1,u_{11})$ be the natural coordinates on the jet space
$J^{2}(2,1)$. Then $({\cal S}_{2})$ is a complex 3-dimensional
submanifold
 in $J^{2}(M \times N)$ defined by the  equation $u_{11} = F(x,u,u_1)$.

{\bf Example 2.} More generally, let $M \subset \cc^n$, $N \subset \cc^m$ be
domains,
$({\cal S})$ be a holomorphic completely overdermined second order system:
$({\cal S}): u^{k}_{x_ix_j} = F^k_{ij}(x,u,u_x)$, $k = 1,...,m$,
$i,j =1,...,n$, $(x,u) \in M \times N$. Denote by $(x,u,u^{k}_i,
u^{k}_{ij})$ the natural coordinates on $J^{2}(n,m)$.
Then $({\cal S}_{2})$ is a complex submanifold in $J^{2}(M \times N)$
defined by
the equations $u^{k}_{ij} = F^{k}_{ij}(x,u,u^{(1)})$ where $u^{(1)} =
(u^k_i)$.

 Since $\pi_M: J^r(M \times N) \longrightarrow M$
also is
a fiber bundle over $M$, every holomorphic map $u: M \longrightarrow N$
defines a
section of this bundle  by $p \longrightarrow j^r_p(u)$.
So $u$ is a holomorphic solution of the system $({\cal S})$ if and only
if the section $p \mapsto j^r_p(u)$ is contained in the variety
$({\cal S}_{r})$.

If $({\cal S}_{r})$ is a regular submanifold of $J^r(M,N)$,
the system $({\cal S})$ is called of {\it maximal rang}. Thus every system
$({\cal S})$  of
maximal rang can be identified with a complex submanifold of the holomorphic
fiber
bundle $\pi_{M \times N}: J^r(M \times N) \longrightarrow M \times N$ and
its solutions
can be identified with sections of the holomorphic fiber bundle $\pi_M:
J^r(M \times N)
\longrightarrow M$.   As we have seen in the above
examples, completely overdetermined systems always are of maximal rang.

\begin{e-def}
 A system
$({\cal S})$ is called {\it locally regular}, if for every point $P \in
J^r(M,N)$ with
the natural projection $\pi_{M \times N}(P) = (p,q) \in M \times N$ there
exists a solution $u(x)$
of $({\cal S})$ holomorphic near $p$ such that $j^r_p(u) = P$.
\end{e-def}

 A holomorphic function $F$ is called an
{\it invariant function} for a one-parameter LTG with an infinitesimal
generator $X$
 if $F(x^*) \equiv F(x)$. It is easy to see that $F(x^*) = e^{tX}F(x)$; this
implies that $F$ is invariant if and only if $X F(x) = 0$. A complex
subvariety $V = \{ F(x) = 0 \}$, where $F$ is a vectror valued holomorphic
function of maximal rang, is called {\it an invariant variety} for a
one-parameter
$LTG$  if $F(x^*) = 0$ when $F(x) = 0$. Clearly, $V$ is an invariant variety
if and only if
$X F(x) = 0$ for every $x \in V$ that is $X$ is a tangent field to $V$.

The importance of these notions explains by the following simple but
fundamental
statement (see for instance \cite{Ol,Po}):

\begin{e-pro}
(The Lie criterion) A vector field $X$ is an infintesimal symmetry of a
locally regular  system $({\cal S})$ of order $r$ and of maximal rang if and
only if the variety  $({\cal S}_{r})$ is   invariant  for the
$r$-prolongation $X^{(r)}$.
\end{e-pro}

It follows by the Cauchy existence theorem that every system of ordinary 
differential equqtions (solved with respect to the highest order derivatives)
is locally regular. In the case of several independent variables we need 
the Frobenius theorem which imposes integrability conditions.

A holomorphic completely overdermined second order PDE system
of the form

\begin{eqnarray*}
({\cal S}): u_{x_ix_j}^k = F_{ij}^k(x,u,u_x), i,j = 1,...,n, k =
1,...,m
\end{eqnarray*}
is always of maximal rang, but in general it is  locally regular. So we need
to assume that it satisfies the integrability condition in the following
sense:
the distribution on the tangent bundle of the jet space $J^1(n,m)$
defined by the differential forms

\begin{eqnarray*}
\omega^k_i = du_i^k - \sum_j F^k_{ij}(x,u,u^{(1)})dx_j, \phi^k = du^k - \sum_iu_i^kdx_i
\end{eqnarray*}
is completely integrable. We call such systems {\it completely integrable}
or
{\it involutive}.  It follows by the Frobenius theorem that
every involutive system is locally regular. Thus, the last proposition is
applicable for this  class of systems.

In the next section we will see that this proposition 
gives an efficient tool for the computation  of infinitesimal
symmetries of
holomorphic completely overdetermined second order involutive systems with 
additional first order relations.

\section{Segre varieties, holomorphic maps and PDE symmetries}

Denote by $Z = (z,w) \in \cc^n \times \cc^m$ the standard coordinates in
$\cc^{n+m}$.
All our considerations will be purely local, so all neighborhoods, domains
etc.
(which we usually even do not mention) always are supposed to be as small as
we need (the most rigorous way is to use the language of germs; following 
the classical tradition we do not employ it in order to avoid useless 
formalizations).
By a {\it real analytic submanifold} ${\cal M}$ of codimension $m$
in $\cc^{n+m}$ we mean
the zero set ${\cal M} = \{ Z: r(Z,\overline{Z}) = 0 \}$ of a
real analytic $\R^m$ - valued map $r = (r^1,...,r^m)$ of maximal rank. Such
a
manifold is called {\it generic} if
$\partial r_1 \wedge ... \wedge \partial r_m \neq 0$. In this paper
we consider generic manifolds only. The {\it holomorphic tangent space}
$H_p({\cal M})$
at a point $p \in {\cal M}$ is the maximal complex subspace of the tangent
space of ${\cal M}$ at $p$. ${\cal M}$ is called Levi nondegenerate at $p$
if the two following conditions hold:
\begin{itemize}
\item[(i)] there exists a linear combination of the Levi forms 
 ${\cal L}_p^j(u,v) = \sum_{j,k} r^i_{z_j
\overline{z}_k}(p)u_jv_k$, $u, v \in H_p({\cal M})$ which is 
a nondegenerate hermitian form on $H_p({\cal M})$ 
\item[(ii)] the forms ${\cal L}^j_p(u,v)$ are $\cc$-linearly
independent.
\end{itemize}
We say that ${\cal M}$ is {\it Levi nondegenerate} if it is
Levi nondegenerate at every point. Often some authors call ${\cal M}$ 
Levi-nondegenerate if a slightly weaker condition holds instead of (i):
the Levi form of ${\cal M}$ (considered as a vector valued hermitian form) 
has the trivial kernel. Our methods can be easily carried to this case (and 
even to a much more general situation). In the present paper we restrict 
ourselves by the consideration of the above class of varieties in order 
to avoid supplementary computations and complications of the notations.

 A map $f:
{\cal M} \longrightarrow {\cal M}$ defined and biholomorphic in
a neighborhood of ${\cal M}$
is called a {\it biholomorphism} or a (biholomorphic) {\it automorphism}
 of ${\cal M}$.
These maps form a group with respect to the composition which is called {\it
the group
of biholomorphisms} or the {\it automorphism group} of ${\cal M}$ and is
denoted by
$Aut({\cal M})$.

The study of automorphism groups of  real submanifolds in $\cc^{n+m}$
is a traditional problem
of the geometric complex analysis and the complex differential geometry.
An important fact here
is that such a group (in the Levi nondegenerate case) is always
a real finite dimensional Lie group.
This phenomenon is due to the intrensic geometry
of a real submanifold  induced
by the complex structure of the ambient space. It has been studied in the
foundator works of E.Cartan \cite{Ca}, N.Tanaka \cite{Ta}, S.S.Chern -
J.Moser
\cite{CM} for the case of real hypesurfaces.
Cartan, Tanaka and Chern study the equivalence problem for a $G$-structure
corresponding to the natural {\cal CR}-structure sitting on a real
hypersurface
in $\cc^{n+1}$ and solve the equivalence problem for this structure using
Cartan's equivalence method for general $G$-structures. In particular, this
gives
a complete list of biholomorphic invariants of a hypersurface. Moser solves
the
equivalence problem via his theory  of a normal form of a real analytic
hypersurface
with respect to the action of local biholomorphisms. This theory gives many
additional
useful information about biholomorphic maps of real hypersurfaces. In
particular, it
leads to an explicit parametrization of the automorphism group. The
approaches of Cartan - Tanaka -Chern and Moser have been developed for
the case of submanifolds of higher codimension in the works of
V.Beloshapka \cite{Be1}, A.Loboda \cite{Lo}, V.Ezhov - A.Isaev - G.Schmalz
\cite{EIS} and other authors.

Another natural approach is to study the Lie algebra $Lie({\cal M})$ of the
automorphism
group of a real analytic manifold ${\cal M}$.
Vector fields in $Lie({\cal M})$    are called {\it infinitesimal
automorphisms} of ${\cal M}$.  The
knowledge of
the Lie algebra allows to refind a neighborhood of the identity in the
automorphism group  via the exponential map i.e. essentially to describe
completely
the group in the local situation. The results in this direction have been
 obtained by   E.Bedford - S.Pinchuk \cite{BP2}, V.Beloshapka
 \cite{Be2}, N.Tanaka \cite{Ta}, A.Tumanov \cite{Tu}, N.Stanton \cite{St},
 and  other authors.

The common feature of all  these works is a direct study of a mixed
"real-complex"
structure of a hypersurface embedded to $\cc^{n+1}$. This leads to
 computations with power series contaning "mixed" terms of the type
$Z^k\overline{Z}^l$ in order to find biholomorphic invariants.  There
is another way to find biholomorphic invariants of a real analytic
submanifold
in $\cc^{n+m}$. For a fixed point $\zeta \in \cc^{n+m}$ close enough to
${\cal M}$ consider the
{\it complex} submanifold $Q(\zeta) = \{ Z: r(Z,\overline\zeta) = 0 \}$.
It is called
the {\it Segre variety} for B.Segre who introduced these objects \cite{Se}.
The basic property
of the Segre varieties is their biholomorphic invariance: for every
automorphism
$f \in Aut({\cal M})$ and any $\zeta$ one has $f(Q(\zeta)) = Q(f(\zeta))$.
 For the approach developed in the present paper,  the utilisation of the
complex conjugation in
the definition of the Segre surface is technically incovenient. So we
consider the
complex  hypersurface $Q^*(\zeta) = Q(\overline\zeta)$. Then for
every $f\in Aut({\cal M})$ one has $f(Q^{*}(\zeta)) = Q^{*}(\overline f
(\overline{\zeta}))$.
Thus, $f$ {\it maps any element of the  family $\{ Q^*(\zeta) \}_{\zeta}$ to
another
one}. This property is crucial for our paper since it can be viewed from the
geometric
PDE point of view. Of course, we still call $Q^*(\zeta)$ the Segre variety
and omit
the star.

The Segre varieties were reintroduced to the modern theory in
the  important
works of  S.S.Chern
\cite{Ch} and S.Webster
\cite{We1} and turned out to be a very useful tool for a study of holomorphic
maps. The theory of Segre varieties has been applied to the study of
 analytic and algebraic 
extension of holomorphic maps  by  M.S.Baouendi -
P.Ebenfelt - L.P.Rothschild \cite{BER1},  K.Diederich - S.Webster
\cite{DW},
 K.Diederich -J.E.Fornaess \cite{DF}, K.Diederich - S.Pinchuk \cite{DP},  S.Webster \cite{We2}. 
 J.Faran \cite{Fa} and S.Webster
 \cite{We3,We4}  also studied related geometric 
invariants.

 M.S.Baouendi - P.Ebenfelt - L.P.Rothschild \cite{BER2, BER3} 
and D.Zaitsev \cite{Za}  
  used the Segre varieties geometry in order to obtain   results
concerning  estimates of dimension and parametrization of automorphism
groups for  various classes of higher codimensional manifolds.

Our approach also makes use of the Segre varieties but the important
difference is that we consider the subject from a more general PDE point of view. It is necessary to stress that the basic idea goes back to the foundators works of B.Segre, E.Cartan and S.Lie's school.

 B.Segre \cite{Se} observed that  in $\cc^2$ the set of Segre varieties of
a Levi nondegenerate real analytic hypersurface ${\cal M}$
(which is called {\it the Segre family} of ${\cal M}$) is a
regular two parameter family of holomorphic curves and so represents the
trajectories
of solutions of a holomorphic second order ordinary differential equation.
The invariance of the Segre family with respect to $Aut({\cal M})$ means
that every biholomorphism of ${\cal M}$ can be considered as a ${\it
symmetry}$
of the differential equation defining its Segre family.

                                       Segre's observation is of fundamental
importance since it links the CR geometry with the PDE theory.

The study of symmetries of a second order ordinary differential
 equation (in some sense,  completed) has been proceeded by S.Lie and his student A.Tresse
\cite{Tr}
(see also \cite{Di}, \cite{Ov}, \cite{Ol2}). In particular, such group is
always a complex Lie group of dimension $\leq 8$; this important 
fact  allowed to B.Segre to
conclude that $Aut({\cal M})$ is a real dimensional Lie group.

The idea of Segre  can be naturally generalized to higher dimension as follows.

First of all, we consider the case where ${\cal M}$ is a real analytic
Levi nondegenerate hypersurface in $\cc^{n+1}$ through the origin.

After a biholomorphic change of coordinates in a neighborhood of the origin
by
the equation $\{ w + \overline{w} + \sum_{j=1}^n \varepsilon_j
z_j\overline{z}_j
+ R(Z,\overline{Z}) = 0 \}$ where $\varepsilon_j = 1$ or $-1$ and $R =
o(\vert Z \vert^2)$.
For every point $\zeta = (\eta_1,...,\eta_n,\omega)$ the corresponding Segre
variety
is given by $w + \omega + \sum_{j=1}^n \varepsilon_jz_j\eta_j + R(Z,\zeta) =
0$.
If we  consider the variables $x_j =z_j$ as independent ones
and the variable $w = u(x)$ as dependent one this equation can be rewritten
in the
form

\begin{eqnarray}
\label{7.1}
 u + \omega + \sum_{j=1}^n \varepsilon_jx_j\eta_j + R(x,\zeta) = 0
\end{eqnarray}
(after an application of the implicit function theorem in 
order to remove $u$ from $R$).
 Taking the derivatives in $x_k$ we obtain the equations

\begin{eqnarray}
\label{7.2}
 u_{x_k} + \varepsilon_k\eta_k + R_{x_k}(x,\zeta) 
= 0,
k = 1,...,n
\end{eqnarray}
The equations (\ref{7.1}), (\ref{7.2}) and the implicit function theorem
imply that $\zeta = \zeta(x,u,u_{x_1},...,u_{x_n})$ is
a holomorphic function; taking again the partial derivatives in $x_j$ in
(\ref{7.2}),  we obtain the following completely overdermined second order holomorphic PDE  system:

\begin{eqnarray*}
({\cal S}_{\cal M}): u_{x_jx_k} = F_{jk}(x,u,u_x), j,k = 1,...,n
\end{eqnarray*}
with $u_x = (u_{x_1},...,u_{x_n})$.  It is very important to point out that
this
system necessarily satisfies the integrability condition of the Frobenius
type.
More precisely, with such a system one can associate the differential forms

\begin{eqnarray*}
\omega_i = du_i - \sum_j F_{ij}(x,u,u^{(1)})dx_j, \phi = du - \sum_iu_idx_i
\end{eqnarray*}
defined on the jet space $J^1(n,1)$. It follows directly from the
representation (\ref{7.1}) of its integral manifolds
that the distribution defined by these forms on the tangent bundle of
$J^1(n,1)$
is completely integrable and so satisfies the Frobenius condition.
The property of biholomorphic invariance of the Segre varieties means that
any biholomorphism of $\Gamma$
transforms the graph of a solution of $({\cal S}_{\cal M})$ to the graph of
another solution,
i.e. is a Lie symmetry of $({\cal S}_{\cal M})$.
 This
naturally leads to a general consideration of a holomorphic involutive
PDE system of the form

\begin{eqnarray*}
({\cal S}_{\cal M}): u_{x_ix_j}^k = F_{ij}^k(x,u,u_x), i,j = 1,...,n, k =
1,...,m
\end{eqnarray*}

Thus, the study of biholomorphisms of real analytic Levi nondegenerate
hypersurfaces can be reduced to the study of symmetries of holomorphic
involutive PDE systems (with one dependent variable). However, the systems
corresponding to Segre families form a very special subclass between
involutive systems since the coefficients of (\ref{7.1}) satisfy additional
conjugation relations due to the fact that the defining function $r$
is real valued. We point out here that the importance of the study of this 
class of PDE systems has been realized by S.S.Chern \cite{Ch} who solved 
the equivalence problem for this class of systems with one dependent variable 
(see also the work of J.Faran \cite{Fa}).

Now consider the higher codimensional case. First of all, we introduce 
the class of PDE systems which plays the major role in the present paper.

Let ${\cal S}$ be a holomorphic second order PDE system with additional 
first order relations of the form

\begin{eqnarray*}
& &u^1_{x_ix_j} = F_{ij}(x,u,u^1_x), i,j = 1,...,n\\ 
& &u^k_x = G^k(x,u,u^1_x), k = 2,...,m
\end{eqnarray*}

In order to simplify the notations we introduce the dependent variables 
$w := u^1$ and $v = (u^2,...,u^m)$ so $u = (w,v)$. Then our system can be
 rewritten in the form 

\begin{eqnarray*}
w_{x_ix_j} = F_{ij}(x,v,w,w_x), v^k_x = G^k(x,v,w,w_x)
\end{eqnarray*}
where we use the notation $w_x = (w_{x_1},...,w_{x_n})$, 
$v^k_x = (v^k_{x_1},...,v^k_{x_n})$, $G^k = (G^k_1,...,G^k_n)$. We will
also use the notation $v_x = (u^2_x,...,u^m_x)$, $G = (G^2,...,G^m)$.

Consider a complex subvariety $\Gamma$ in the jet space $J^1(n,m)$ 
defined by $(x,u,u^{(1)}): v^{(1)} = G(x,u,w^{(1)})$ in the natural 
coordinates. Then $(x,u,w^{(1)})$ are holomorphic local coordinates 
on $\Gamma$ and we may consider the 1-forms defined on $\Gamma$ as follows:

\begin{eqnarray*}
\omega_i = dw_i - \sum_j F_{ij}(x,u,w^{(1)})dx_j, 
\phi^1 = dw - \sum_j w_jdx_j,
\phi^k = dv^k - \sum_j G_j^k(x,u,w^{(1)})dx_j, k > 1
\end{eqnarray*}

We say that the system $({\cal S})$ is {\it completely integrable}
 or {\it involutive} 
if the distribution defined by these forms on the tangent bundle of $\Gamma$
is completely integrable that is satisfies the Frobenius condition. It 
follows by the Frobenius theorem that if $({\cal S})$ is involutive then 
it is locally regular i.e. for every point  of the complex submanifold 

\begin{eqnarray*}
{\cal S}_2: w_{ij} = F(x,u,w^{(1)}), v^{(1)} = G(x,u,w^{(1)})
\end{eqnarray*}
of $J^2(n,m)$ there exists a solution of $({\cal S})$ whose jet coincides with 
this point. In view of the Frobenius criterion  the graphs of solutions of $({\cal S})$ form a holomorphic foliation of $\Gamma$ with n-dimensional leafs
 and depending on (n+m)-parameters if and only if $({\cal S})$ is involutive.

Let now ${\cal M}$ be a Levi nondegenerate
quadric in $\cc^{n+m}$ given by $w_k + \overline{w}_k =
 <L^k(z),\overline{z}>$, $k = 1,...,m$ where every $L^k$ is a
hermitian operator on $\cc^n$ and $<z,\zeta> = \sum_{j=1}^n
z_j\zeta_j$.  We can assume 
that the hermitian form $<L^1(z),\overline{z}>$ is nondegenerate.
For $(z,\omega) \in \cc^n \times \cc^m$ the corresponding Segre variety is 
$Q(\zeta,\omega) = \{ (z,w): w_k + \omega_{k} = <L^k(z),\zeta> \}$.
If we consider $x:= z$ as independent variables and $u: = w$ as dependent,
then
$Q(\zeta,\omega)$ is a graph of $u$: $Q(\zeta,\omega) = \{ (x,u): u^k  + \omega_k = <L^k(x),\zeta> \}$.

Let us construct a PDE system with a general solution given by the above
family.
First of all, clearly we have the equations $u^k_{x_ix_j} = 0$, for every
$k,i,j$. However, in general this is not enough since our family of
solutions
depends only on $n +m$ parameters and so we need to look for another
relations.
Considering the first partial derivatives we obtain the following system
of linear algebraic equations for $\zeta$: $u^1_{x_i} = <L^1(x)_{x_i},\zeta>$.
Since  the rank of this system is equal to $n$, we get $\zeta = N u^1_x$, where
 $N$ is an $n \times n$
matrix. Set as above  $v=(u^2,...,u^m)$ and $w = u^1$.  We obtain that 
$u^k_x = A^k u^1_x $, $k = 2,...,m$, 
where every $A^k$ is a matrix. Therefore, we obtain the following PDE system:

\begin{eqnarray*}
u^1_{x_ix_j} = 0, u^k_x = A^k u^1_x, k = 2,...,m
\end{eqnarray*}
whose sets of solutions coincides with the Segre family of ${\cal M}$.

This construction can be immediately generalized to any Levi nondegenerate
real analytic submanifold. Indeed, let ${\cal M}$ be a real analytic
Levi nondegenerate submanifold in $\cc^{n +m}$ through the origin. Then
in a neighborhood of the origin it can be represented in the form 
$w_k + \overline{w}_k = <L^k(z),\overline{z}> + o(\vert Z\vert^2)$, 
$k = 1,...,m$.
For $(z,\omega) \in \cc^n \times \cc^m$ the corresponding Segre variety is
$Q(\zeta,\omega) = \{ (z,w): w_k + \omega_{k} = <L^k(z),\zeta> +
R^k(z,\zeta,\omega)\}$
where $R^k$ contains no term of order $\leq 2$ (after an application 
of the implicit function theorem if it is necessary).
Consider $x:= z$ as independent variables and $u: = w$ as dependent, then
$Q(\zeta,\omega)$ is a graph of $u$:

\begin{eqnarray}
\label{7.3}
Q(\zeta,\omega) = \{ (x,u): u^k  + \omega_k = <L^k(x),\zeta> +
R^k(x,\zeta,\omega)\}
\end{eqnarray}

Considering the first partial derivatives we obtain the following system :

\begin{eqnarray}
\label{7.4}
u^k_{x_i} = <L^k(x)_{x_i},\zeta> + R^k_{x_i}(x,\zeta,\omega)
\end{eqnarray}
Applying the implicit function theorem to (\ref{7.3}), (\ref{7.4}) we get
that $(\zeta,\omega) = \varphi(x,u,w_x)$, where $\varphi$ is a holomorphic
function. It is worth to point out that the implicit function theorem allows
to compute by recursion a term of any order in the expansion of $\varphi$,
so our method is totally constructive. Using $\varphi$ in order to exclude
the parameters $\zeta$, $\omega$ from those equations of (\ref{7.4}) which
are not
used yet, we obtain holomorphic equations of the form $u^k_x = A^k u^1_x +
\psi(x,u,u^1_x)$, $k = 2,...,m$
with holomorphic function $\psi$ without terms of order $\leq 1$.

Next, we consider the second order partial derivatives $u^1_{x_ix_j} = 
R^1_{x_ix_j}(x,\zeta,\omega)$
and replace  $\zeta$, $\omega$ by $\varphi$. 
 We obtain the holomorphic equations $w_{x_ix_j} = F_{ij}(x,u,w_x)$.
Thus,
 finally we obtain that $u(x)$ satisfy the following holomorphic PDE system:

 \begin{eqnarray*}
& &({\cal S}_{\cal M}): w_{x_ix_j} = F_{ij}(x,u,u_x), i \leq j, v^k_x = A^k w_x
+ G^k(x,u,w_x)
\end{eqnarray*}

Since the solutions of this system (given by (\ref{7.3})) depend on $(n+m)$ parameter, it follows by the Frobenius theorem that this system is involutive 
(in particular, (\ref{7.3}) represents all solutions of this system).

 The
biholomorphic invariance of the Segre family of ${\cal M}$ means that every
biholomorphism of ${\cal M}$ is a symmetry of the constructed PDE system.

Therefore, in the case where $Sym({\cal S}_{\cal M})$ is a finite
dimensional
Lie group, $Aut({\cal M})$ is its finite dimensional real Lie subgroup (since 
it is obviouisly closed).
In order to obtain a precise  estimate of its dimension,
we recall the following useful observation due to
E.Cartan \cite{Ca}. Let a holomorphic vector field $X$ be an infinitesimal
generator of $Aut({\cal M})$ (this means that we consider the real time $t$
in the corresponing Lie series). This is
equivalent to the fact that $Re X$ is a tangent vector field to ${\cal
M}$. 
On the other hand,
$X$ is an infinitesimal symmetry of $({\cal S}_{\cal M})$. Indeed, every
biholomorphism
from the corrseponding real one-parameter group takes an element of the Segre
family to another one, so $X$ is tangent to $Sym({\cal S}_{\cal M})$
considered
as a real Lie group; but since $X$ is a holomorphic vector field, it
is 
necessarily in $Lie({\cal S}_{\cal M})$.
It is clear that if ${\cal M}$ is Levi nondegenerate, the field $Re(iX)$
cannot be
tangent to ${\cal M}$ simultaneously with $Re X$ i.e.  $Lie({\cal M})$
is a totally real subspace of $Lie({\cal S}_{\cal M})$. Therefore, the real
dimension of $Aut({\cal M})$ is majorated by  the complex dimension
of $Lie({\cal S}_{\cal M})$.

We stress again that  quite similarly to the hypersurface case, systems
defining the Segre families form a very special subclass of the class of
second order holomorphic involutive systems with first order relations.

We have proved the following

\begin{e-pro}
 The Segre family of a real analytic Levi nondegenerate submanifold
 ${\cal M}$
of $\cc^{n+m}$ is a general solution of a holomorphic second order
completely overdetermined
involutive PDE system with $n$ independent and one dependent variables and
first order relations. This system
is canonically associated with ${\cal M}$ and is denoted by $({\cal S}_{\cal
M})$.

 If $Sym({\cal S}_{\cal M})$ is a finite dimensional
complex Lie group, then $Aut({\cal M})$ is its real Lie subgroup
embedded to $Sym({\cal S}_{\cal M})$ as a totally real submanifold.
\end{e-pro}

We conclude this section by some examples. It is easy to show (see 
\cite{Lo}) that every 6-dimensional quadric in $\cc^4$ is linearly 
equivalent to one of the following quadrics:

\begin{eqnarray*}
& &{\cal M}^1: w_1 + \overline{w} = z_1\overline{z_1} + z_2\overline{z}_2,
w_2 + \overline{w}_2 = z_1\overline{z}_1 - z_2\overline{z}_2,\\
& &{\cal M}^1: w_1 + \overline{w}_1 = z_1\overline{z}_1 - z_2\overline{z}_2,
w_2 + \overline{w}_2 = z_1\overline{z}_2 + z_2\overline{z}_1,\\
& &{\cal M}^3: w_1 + \overline{w}_1 = z_1\overline{z}_2 + z_2\overline{z}_1,
w_2 + \overline{w}_2 = z_1\overline{z}_1
\end{eqnarray*}

Considering independent variables $x = z$ and dependent variables 
$u = w$ we get that the systems defining the corresponding Segre families 
are

\begin{eqnarray*}
& &({\cal S})^1: u^1_{x_ix_j} = 0, i,j = 1,2, u^2_{x_1} = u^1_{x_1},
 u^2_{x_2} = - u^1_{x_2},\\
& &({\cal S})^2: u^1_{x_ix_j} = 0, i,j = 1,2, u^2_{x_1} =- u^1_{x_2}, 
u^2_{x_2} =  u^1_{x_1},\\
& &({\cal S})^3: u^1_{x_ix_j} = 0, i,j = 1,2, u^2_{x_1} = u^1_{x_2},
 u^2_{x_2} = 0
\end{eqnarray*}

In the next two sections we develop a general approach in order to study
infinitesimal symmetries of  
second order holomorphic involutive systems with first order relations. 
Much more advanced tools can be found in \cite{Po,Po2}; we only adapt for our 
case a very elementary part of the general theory.

\section{Completely integrable systems, their deformations and
infinitesimal symmetries }

Consider a holomorphic second order involutive  PDE system $({\cal S}_0)$

\begin{eqnarray*}
& &({\cal S}_0):
w_{x_{i_1}x_{i_2}} = F_{i_1i_2}( x,u, w_x),
i_1,i_2 = 1,...,n, \mu = 1,...,m,\\
& &v^k_{x_j} = G^k_j(x,u,w_x), k = 2,...,m, j = 1,...,n
\end{eqnarray*}
with $n$ independent variables $x$ and $m$ dependent variables $u = (w,v) \in \cc \times \cc^{m-1}$.

By a completely integrable holomorphic deformation  of the system $({\cal
S}_0)$
 we mean a PDE system of the form

\begin{eqnarray*}
({\cal S}^{\varepsilon}): w_{x_ix_j} = F_{i_1i_2}(\varepsilon, x,u, w_x),
v_x =  G(\varepsilon,x,u,w_x)
\end{eqnarray*}
where $F^{\mu}_{i_1i_2}$, $G$ are holomorphic functions in $(x,u,w_x)$
and real analytic with respect to a (vectorvalued) parameter
$\varepsilon$; they  satisfy
$F_{i_1i_2} \vert \{\varepsilon = 0 \} \equiv F_{ij}$,
$G \vert \{ \varepsilon = 0\} = G$ and are such that this system is completely
integrable for every fixed  $\varepsilon$.

For every $\varepsilon$ we can consider all first order partial 
derivatives of the equations $v_x = G(\varepsilon,x,u,w_x)$ and then 
substitute $w_{x_ix_j} = F_{ij}(x,u,w_x)$ in order to remove the second 
order derivatives of $w$ in the right sides. The obtained PDE system has
the form 

\begin{eqnarray*}
& &({\cal S}^{\varepsilon}): u^k_{x_ix_j} =
 F^k_{i_1i_2}(\varepsilon, x,u, w_x), k = 1,...,m, i_1,i_2 = 1,...,n,\\
& &v_x =  G(\varepsilon,x,u,w_x)
\end{eqnarray*}
and obviously has the same space of solutions as the initial system, 
so has the same symmetry group. We will work with this system.

In order to study $Lie({\cal S}^{\varepsilon})$ we apply the general Lie
method to the deformed system
$({\cal S}^{\varepsilon})$. This system defines a complex
subvariety $({\cal S}^{\varepsilon}_{2})$ of the jet space $J^2(n,m)$
given by the equations

\begin{eqnarray*}
 u^{\mu}_{i_1i_2} = F^{\mu}_{i_1i_2}(\varepsilon,x,u,u^{(1)}), v^{(1)} = G(\varepsilon,x,u,w^{(1)})
\end{eqnarray*}
and in view of the integrability condition this system is locally regular.
Therefore the Lie criterion implies that
$X = \sum \theta^j \frac{\partial}{\partial x_j} + \sum \eta^{\mu}
\frac{\partial}{\partial u^{\mu}}$ is in $ Lie({\cal S}^{\varepsilon})$ if and only if $X^{(2)}$ is tangent
to $({\cal S}^{\varepsilon}_{2})$. This is equivalent to the
following
equations:

\begin{eqnarray*}
& & X^{(2)} u^{\mu} = X^{(2)} (F^{\mu}_{i_1i_2}(\varepsilon,x,u,u^{(1)})) =
X^{(1)}(F^{\mu}_{i_1i_2}(\varepsilon,x,u,u^{(1)}))\\
& & X^{(2)}(v^{(1)} - G(\varepsilon,x,u,w^{(1)})) = 0, (x,u,u^{(1)},u^{(2)}) \in ({\cal S}^{\varepsilon}_{2})
\end{eqnarray*}

Clearly, {\it this is a linear condition on the coefficients $\theta$,
$\eta$
of $X$ and their partial derivatives up to the second order}. We explain now
how to construct explicitely the corresponding linear second order
PDE system with holomorphic coefficients for $\theta$, $\eta$ 
 equivalent to this condition.

Set $\hat{\eta}^{\mu}_{i_1i_2} = \eta^{\mu}_{i_1i_2} -
\Lambda^{\mu}_{i_1i_2}$.
Then we have

\begin{eqnarray*}
 \hat{\eta}^{\mu}_{i_1i_2} = - \Lambda^{\mu}_{i_1i_2} +
X^{(1)}(F^{\mu}_{i_1i_2}(\varepsilon,x,u,u^{(1)}),
(x,u,u^{(1)},u^{(2)}) \in ({\cal S}_2)^{\varepsilon}
\end{eqnarray*}

Set $L_2 = \{ (x,u,u^{(1)},u^{(2)}): u^{\mu}_{i_1i_2} =
F^{\mu}_{i_1i_2}(\varepsilon,x,u,u^{(1)})
\}$ and $L_1 =
\{ (x,u,u^{(1)},u^{(2)}): v^{(1)}= G(\varepsilon,x,u,w^{(1)}) \}$,
so $({\cal S}^{\varepsilon}_{2}) = L_1 \cap L_2$.

Using the equalities $u^{\mu}_{i_1i_2} =
F^{\mu}_{i_1i_2}(\varepsilon,x,u,u^{(1)})$ we replace
$u^{\mu}_{i_1i_2}$ by $F^{\mu}_{i_1i_2}$ in $\Lambda^{\mu}_{i_1i_2}$
and denote obtained expressions by $\hat\Lambda^{\mu}{i_1i_2}$. We point out
that
they are linear in $\partial \theta$, $\partial \eta$ (the vector functions formed by all first order partial derivatives of $\theta^j$, $\eta^{\mu}$). 
We get the equations

\begin{eqnarray}
\label{4.1}
\hat\eta^{\mu}_{i_1i_2} \vert L_2=
- \hat\Lambda^{\mu}_{i_1i_2}(\varepsilon,x,u,u^{(1)}, \partial \theta,
\partial \eta) 
+ \phi^{\mu}_{i_1i_2}
(\varepsilon,x,u,u^{(1)},\theta,\eta,\partial \theta,\partial \eta)
\end{eqnarray}
where holomorphic functions $ \phi^{\mu}_{i_1i_2}(\varepsilon,x,u,u^{(1)},\theta,\eta,\partial
 \theta,\partial \eta) = X^{(1)} F^{\mu}_{i_1i_2}(\varepsilon,x,u,u^{(1)})$
are linear with respect to $\theta$, $\eta$, $\partial \theta$, $\partial
\eta$.
On the other hand, $\hat\eta^{\mu}_{i_1i_2}\vert L_2 = \sum_{\vert \alpha \vert \leq 3}A^{\mu}_{i_1i_2\alpha}[u^{(1)}]^{\alpha}$
where the coefficients $A^{\mu}_{i_1i_2\alpha}$ are integer
linear combinations with constant coefficients of second order
partial derivatives of $\theta$, $\eta$ (of course, we suppose that
$A^{\mu}_{i_1i_2\alpha}$ are defined for every $\alpha$ allowing them to
vanish
identically).

Next we need to restrict  our expressions  on $L_1$: 
$\hat\eta^{\mu}_{i_1i_2} \vert ({\cal S}^{\varepsilon})_{2} =
\sum_{\vert \beta \vert }
B^{\mu}_{i_1i_2\beta}[w^{(1)}]^{\beta}$
where $B^{\mu}_{i_1i_2\beta} = \sum_{\alpha} b^{\mu\alpha}_{i_1i_2\beta}
 A^{\mu}_{i_1i_2\alpha}$
and the coefficients $b^{\mu\alpha}_{i_1i_2\beta}$ are holomorphic functions
in $(\varepsilon, x, u)$. Therefore, every $B^{\mu}_{i_1i_2\beta}$ is a
linear
combination of the second order partial derivatives of $\theta$,
$\eta$ of 
the form $B^{\mu}_{i_1i_2\beta} = \sum_{j; \vert \gamma \vert = 2} c^j_{A\gamma}\partial^{\gamma} \theta_j + \sum_{k; \vert \delta \vert = 2}
d^{k}_{A\delta}\partial^{\delta} \eta^k$
where we write $ A = (\mu,i_1,i_2,\beta)$ for simplicity of notations and
the
coefficients are holomorphic functions in $(\varepsilon, x, u)$.

Developing the right sides of (\ref{4.1}) into power series with respect to
$u^{(1)}$ we obtain the series of the form $\sum_{\alpha}
f^{\mu}_{i_1i_2\alpha}(\varepsilon,x,u,\theta,\eta,\partial \theta,
\partial \eta) [u^{(1)}]^{\alpha}$
where the holomorphic coefficients
$f^{\mu}_{i_1i_2\alpha}(\varepsilon,x,u,\theta,\eta,\partial \theta,
\partial \eta)$ are linear with respect to $\theta$, $\eta$, $\partial
\theta$,
$\partial \eta$.

Replacing here $v^{(1)}$ by $G(\varepsilon,x,u,w^{(1)})$  and developing the
obtained expressions in power series in $w^{(1)}$, we obtain that $(\ref{4.1})$
implies $\sum_{\vert \beta \vert } B^{\mu}_{i_1i_2\beta}[w^{(1)}]^{\beta} =
\sum_{\beta} p^{\mu}_{i_1i_2\beta}(\varepsilon,x,u,\theta,\eta,\partial
\theta,
\partial \eta) [w^{(1)}]^{\beta}$
 that is $B^{\mu}_{i_1i_2\beta}  =
p^{\mu}_{i_1i_2\beta}(\varepsilon,x,u,\theta,\eta,\partial \theta,
\partial \eta)$
for any $\mu$, $i_1 \leq i_2$, $\vert \beta \vert $ where the right sides are
linear
with respect to $\theta$, $\eta$, $\partial \theta$, $\partial\eta$.
 By the Noetherian property, there exists a finite number $N$
(independent of $\varepsilon$) such that this
equivalent
 to $B^{\mu}_{i_1i_2\beta}  =
p^{\mu}_{i_1i_2\beta}(\varepsilon,x,u,\theta,\eta,\partial \theta,
\partial \eta), \vert \beta \vert \leq N$

We get a linear PDE system of the form  (using the notation
$z = (x,u)$):

\begin{eqnarray}
\label{4A}
\sum_{\vert \alpha \vert =2}
(a_{j\alpha}^{t}(\varepsilon,z)\partial^{\alpha}
 \theta_j + \sum_{\vert \beta \vert = 2}b^{t}_{k\beta}(\varepsilon,z)
\partial^{\beta}\eta^{k} = c^t(\varepsilon,z,\theta, \eta,\partial \theta,
\partial \eta), t = 1,...,N_1
\end{eqnarray}

\begin{eqnarray}
\label{4A1}
\sum_{j,k}d^t_{j,k}(\varepsilon,z)\frac{\partial \theta_j}{\partial z_k}
+
\sum_{i,l} e^t_{i,l}(\varepsilon,z)\frac{\partial \eta^s}{\partial z_l} =
f^t(\varepsilon,z,\theta,\eta), t= 1,...,N_2
\end{eqnarray}

\begin{eqnarray}
\label{4A2}
\sum_p g^t_p(\varepsilon,z)\theta_p + \sum_q h^t_q(\varepsilon,z)\eta^q =
l^t(\varepsilon,z), t = 1,...,N_3
\end{eqnarray}
where the right sides are linear functions in $\theta$, $\eta$,
$\partial \theta$, $\partial \eta$ (recall that our initial tangency conditions
are linear with respect to $\theta$, $\eta$. Therefore, the right sides
do not contain terms without $\theta$, $\eta$ and their derivatives; in particular, $l^t$ vanishes identically).

Now we proceed quite similarly with the equations

\begin{eqnarray}
\label{4.4}
X^{(1)}v^{(1)} = X^{(1)}G(\varepsilon,x,u,w^{(1)})
\end{eqnarray}

We have $\eta^{\mu}_{i} = \sum_{\vert \alpha \vert \leq 2}
Q^{\mu}_{i\alpha}[u^{(1)}]^{\alpha}$
where $Q^{\mu}_{i\alpha}$ are linear combinations of second order partial
derivatives of $\theta$, $\eta$ with constant coefficients.  The
equations $(\ref{4.4})$ can be rewritten in the form

\begin{eqnarray*}
& &\eta^{\mu}_{i} = X^{(1)} G^{\mu}_i(\varepsilon,x,u,w^{(1)}) =
\sum_{p=1}^n \eta_p^1\psi^{\mu}_{ip}(\varepsilon,x,u,w^{(1)}) +
\phi^{\mu}_{i}(\varepsilon,x,u,w^{(1)},\theta,\eta),\\
& &\mu =2,...,m, i =1,...,n
\end{eqnarray*}
where

\begin{eqnarray*}
& & \psi^{\mu}_{ip}(\varepsilon,x,u,w^{(1)}) =
\frac{\partial G^{\mu}_i}{\partial u_p^1} =
\sum_{\delta}\Psi_{\delta}[w^{(1)}]^{\delta},\\
& & \phi^{\mu}_{i}(\varepsilon,x,u,w^{(1)},\theta,\eta) = \sum_j \theta_j
\frac{\partial G^{\mu}_i}{\partial x_j} + \sum_k \eta^k
\frac{\partial G^{\mu}_i}{\partial u^k} =
\sum_{\delta}\Phi_{\delta}[w^{(1)}]^{\delta}
\end{eqnarray*}

In particular, the functions
$\phi^{\mu}_{i}(\varepsilon,x,u,w^{(1)},\theta,\eta)$ are
linear with respect to $\theta$, $\eta$.

This is equivalent to the equalities

\begin{eqnarray}
\label{4.5}
\sum_{\vert \alpha \vert \leq 2} Q^{\mu}_{i\alpha}[u^{(1)}]^{\alpha} =
\sum_{\vert \alpha \vert \leq 2; p=1,...,n}
Q_{p\alpha}^1\psi^{\mu}_{ip}(\varepsilon,x,u,w^{(1)})[u^{(1)}]^{\alpha}
 + \phi^{\mu}_{i}(\varepsilon,x,u,w^{(1)},\theta,\eta)
 \end{eqnarray}
under the condition $v^{(1)} = G(\varepsilon,x,u,w^{(1)}) =
\sum_{\gamma}g_{\gamma}(\varepsilon,x,u)[w^{(1)}]^{\gamma}$.

Substituting these power series into $(\ref{4.5})$ we get the following
equality (using the vector notation): 
$\sum_{\beta} T_{\beta}[w^{(1)}]^{\beta} = \sum_{\beta} S_{\beta}
[w^{(1)}]^{\beta} +
\sum_{\beta} P_{\beta} [w^{(1)}]^{\beta}$
of power series with vector valued coefficients $T_{\beta}$, $S_{\beta}$
which are linear
combinations of first order partial derivatives of $\theta$, $\eta$ with
coefficients holomorphic in $(\varepsilon,x,u)$ and
$P_{\beta}(\varepsilon,x,u,\theta,\eta)$ being linear in $\theta$, $\eta$.
So we obtain
the following system of the equations: $T_{\beta} - S_{\beta} -
P_{\beta} = 0$ which in view of the Noetherian condition
is equivalent to $T_{\beta} - S_{\beta} - P_{\beta} = 0$, $ \vert \beta
\vert \leq N$ for a finite $N$.

So we have a first order linear system of equations:

\begin{eqnarray}
\label{4B}
\sum_{j,k}\hat a^t_{j,k}(\varepsilon,z)\frac{\partial \theta_j}{\partial
z_k} +
\sum_{i,l} \hat b^t_{i,l}(\varepsilon,z)\frac{\partial \eta^s}{\partial z_l}
=
\hat c^t(\varepsilon,z,\theta,\eta), t= 1,...,N_4,
\end{eqnarray}
\begin{eqnarray}
\label{4B1}
\sum_p \hat d^t(\varepsilon,z)\theta_p + \sum_q \hat
e^t(\varepsilon,z)\eta^q =
\hat f^t(\varepsilon,), t = 1,...,N_5
\end{eqnarray}

As above, the right sides does not contain terms without $\theta$, $\eta$ (
for instance, $\hat{f} \equiv 0$).

 We have proved the following

\begin{e-theo}
\label{theo4.1}
The vector field $X$ defines an infinitesimal symmetry of $({\cal
S}^{\varepsilon})$
if and only if its coefficients satisfy the united system $(\ref{4A})$,
$(\ref{4A1})$, $(\ref{4A2})$,
$(\ref{4B})$, $(\ref{4B1})$. The Lie algebra $Lie({\cal S}^{\varepsilon})$ is finite
dimensional
if and only if the linear space of holomorphic solutions of this united
system is finite dimensional.
\end{e-theo}

The constructed linear holomorphic PDE system  is called  {\it the (infinitesimal) Lie equations} associated
with $({\cal S}^
{\varepsilon})$.

As an important example, let us construct the Lie equations for a PDE system
of the form

\begin{eqnarray}
\label{3.1}
u^k_{x_ix_j} = 0, i,j = 1,...,n, k = 1,...,m
\end{eqnarray}
\begin{eqnarray}
\label{3.2}
 v^k_x = M^kw_x, k = 2,...,m
\end{eqnarray}

We call such a system by a {\it flat system with relations} $({\cal
S}_{flat})$. Obviously, such a system is involutive.

The variety $({\cal S}_{flat})_{2}$ defined by $(S_{flat})$ is given
by the
equations

\begin{eqnarray*}
& &u^k_{ij} = 0, k = 1,...,m, i,j = 1,...,n\\
& & v^{(1)} = Mw^{(1)}
\end{eqnarray*}
where the matrix $M$ is formed by the matrices $M^k$ as vertical blocks.
Let a
vector field $X = \sum_{j= 1}^n \theta^j \frac{\partial}{\partial x_j} + \sum_{\mu =1}^m
\eta^{\mu}
\frac{\partial}{\partial u^{\mu}}$
be in $Lie({\cal S}_{flat})$ i.e. an infinitesimal symmetry of $({\cal
S}_{flat})$.

Since our system is locally regular and of maximal rank, $X \in Lie({\cal
S})_{flat}$
if and only if $X^{(2)}$ is tangent to $({\cal S}_{flat})_{2}$ i.e.

\begin{eqnarray*}
& &X^{(2)} u^{\mu}_{i_1i_2} = 0, i_1, i_2 = 1,...,n, \mu = 1,...,m\\
& &X^{(2)} (v^{(1)} - Mw^{(1)}) = X^{(1)}(v^{(1)} - Mw^{(1)}) = 0,\\
& &u^{\mu}_{i_1i_2} = 0, v^{(1)} = Mw^{(1)}
\end{eqnarray*}

The first line equations imply that

\begin{eqnarray}
\label{3.3}
\eta^{\mu}_{i_1i_2} = 0, (x,u,u^{(1)},u^{(2)}) \in ({\cal S}_{flat})_{2}
\end{eqnarray}
for any $\mu$ and any $i_1 \leq i_2$. We point out also that the equations
$u^{\mu}_{i_1i_2} = 0$ imply

\begin{eqnarray}
\label{3.4}
\Lambda^{\mu}_{i_1i_2} = 0
\end{eqnarray}
Set $L_2 = \{ (x,u,u^{(1)},u^{(2)}): u^{\mu}_{i_1i_2} = 0\}$ and $L_1 =
\{ (x,u,u^{(1)},u^{(2)}: v^{(1)} = Mw^{(1)} \}$, so $({\cal S}_{flat})_{2} = L_1 \cap
L_2$.

In view of $(\ref{3.4})$

\begin{eqnarray}
\label{3.5}
\eta^{\mu}_{i_1i_2}\vert L_2 = \sum_{\vert \alpha \vert \leq 3}
A^{\mu}_{i_1i_2\alpha}
[u^{(1)}]^{\alpha}
\end{eqnarray}
where the coefficients $A^{\mu}_{i_1i_2\alpha}$ are integer linear
combinations of second order
partial derivatives of $\theta$, $\eta$.

Next we need to restrict the polynomials $(\ref{3.5})$ on $L_1$. Replacing
$v^{(1)}$ by $Mw^{(1)}$
in $(\ref{3.5})$ we obtain $\eta^{\mu}_{i_1i_2} \vert ({\cal S}_{flat})_{2} = \sum_{\vert \beta \vert
\leq 3}
B^{\mu}_{i_1i_2\beta}w^{\beta}$
where $B^{\mu}_{i_1i_2\beta} = \sum_{\alpha} b^{\mu\alpha}_{i_1i_2\beta}
A^{\mu}_{i_1i_2\alpha}$
and the coefficients $b^{\mu\alpha}_{i_1i_2\beta}$ are polynomials of degree
$\leq 3$
of elements of the matrix $M$. Therefore, every $B^{\mu}_{i_1i_2\beta}$ is a
linear
combination of the second order partial derivatives of $\theta$, $\eta$: 
$B^{\mu}_{i_1i_2\beta} = \sum_{j; \vert \gamma \vert = 2} c^j_{A\gamma}
\partial^{\gamma} \theta_j + \sum_{k; \vert \delta \vert = 2}
d^{k}_{A\delta}
\partial^{\delta} \eta^k$
where we write $ A = (\mu,i_1,i_2,\beta)$ for simplicity of notations.

Therefore, the equations $(\ref{3.3})$ are equivalent to 
$B^{\mu}_{i_1i_2\beta} = 0$
for all $\mu,\beta,i_1,i_2$.

Now we proceed quite similarly with the equations $X^{(1)}(v^{(1)} - Mw^{(1)}) = 0$
which
are equivalent to the conditions $\eta^{\mu}_{i}\vert L_1 = 0$, 
$\mu = 2,...,m$, $i = 1,...,n$.

 We may write $\eta^{\mu}_{i}  \vert L_1= \sum_{j, \vert \alpha \vert = 1}
 e^{\mu j}_{i\alpha}\partial^{\alpha}\theta_j +
\sum_{k, \vert \beta \vert = 1} f^{\mu j}_{i\beta}\partial^{\beta}\eta^k$
where the coefficients $e^{\mu j}_{i \alpha}$, $f^{\mu j}_{i \beta}$ are
polynomials in the elements of  $M$ of degree $\leq 2$.

Consider now the following second PDE system $({\cal R}_2)$ for the unknown
vector function $\tau := (\theta,\eta)$:

\begin{eqnarray}
\label{3.11}
B^{\mu}_{i_1i_2\beta}= 0, (\eta^{\mu}_i \vert L_1) = 0
\end{eqnarray}
for all $\mu,i_1 \leq i_2, \beta$. This is a linear
second order PDE system with constant coefficients which represents  the
Lie equations for $({\cal S}_{flat})$. We emphasize the very important 
property of this system: {\it  every equation of second (resp. first) order contains only the second (resp. first) order partial derivatives}.

In the next section we recall some general
properties of  linear PDE systems with holomorphic coefficients useful
for a study of the Lie equations.

\section{Symbols, prolongations and solutions of linear systems}

In this section we adapt general methods of the formal PDE theory for 
our case. Much more general methods and tools can be found in \cite{Po,Po2}.

As usual, by a holomorphic linear PDE system of order $q$ with
$n$ independent variables $y$ and $m$ dependent variables $\tau$ we mean a
system
of the form

\begin{eqnarray*}
({\cal R}_q): \sum_{j = 1,...,m;\vert \alpha \vert \leq q}
a^{i}_{j\alpha}(y) \partial^{\alpha} \tau^j = 0, i= 1,...,s
\end{eqnarray*}
where $a^i_{j\alpha}$ are holomorphic functions. We use the same notation
for the subvariety in the jet space $J^q(n,m)$ corresponding to this system:

\begin{eqnarray*}
({\cal R}_q): \sum_{j = 1,...,m; \vert \alpha \vert \leq q}
a^{i}_{j\alpha}(y) \tau^j_{\alpha} = 0, i= 1,...,s
\end{eqnarray*}

A (holomorphic) solution  of such a system is a  function
$\tau(y)$ holomorphic on a domain $D$ of definition of the coefficients such
that
$j^q_x(\tau) \in ({\cal R}_q)$ for every $x \in D$. We denote by $Sol({\cal
R}_q)$
the vector space of the solutions of $({\cal R}_q)$.

{\it The symbol} $G_q(y^0)$ of $({\cal R}_q)$ at a point $y^0$
is a linear subspace of the complex affine space with coordinates
$v^j_{\alpha}$, $j = 1,...,m$, $\vert \alpha \vert = q$, $\alpha_1 \leq ...\leq \alpha_q$ $\alpha_i \in 
\{ 1,...,n \}$,  defined by

\begin{eqnarray*}
(G_q): \sum_{j = 1,...,m; \vert \alpha \vert = q}
a^{i}_{j\alpha}(y^0) v^j_{\alpha} = 0, i= 1,...,s
\end{eqnarray*}

The $r$-prolongation $({\cal R}_{q +r})$ of  $({\cal R}_q)$ is a linear
system which we get  if we add to $({\cal R}_q)$ the equations
obtained by taking  all the partial derivatives of order $\leq r$
in every equation of  $({\cal R}_q)$ , that is

\begin{eqnarray*}
({\cal R}_{q+r}): \sum_{j = 1,...,m;  \vert \alpha \vert \leq q}
\partial^{\beta}(a^{i}_{j\alpha}(y) \partial^{\alpha} \tau^j) = 0, i= 1,...,s, \vert \beta \vert \leq r
\end{eqnarray*}

Obviously, it has the same space of solutions. The symbol of  $({\cal R}_{q +r})$ is denoted by $G_{q +r}(y^0)$.

The system $({\cal R}_{q})$  is called of   {\it finite type} at $y^0$ if
$G_{q +r}(y^0) = \{ 0 \}$ for some $r$. If a system is
of finite type  at every point, we say simply that it is of finite type.
The smallest $r$ with this property is called the type of    $({\cal
R}_{q})$ and is denoted
by $type({\cal R}_q)$.

\begin{e-theo}
Suppose that $({\cal R}_{q})$ is of finite type at some point $y^0$.
Then the dimension of the space of solutions of $({\cal R}_{q})$ 
holomorphic in a neighborhood of $y^0$  is finite.
\end{e-theo}

\proof The fact that $G_{q +r}(y^0) = \{ 0 \}$ for some $r$ implies that
$({\cal R}_{q+r})$ contains a subsystem which can be solved 
with respect to all partial derivatives of order $q+r$ and so can be represented in the form (in a neighborhood of
$y^0$):

\begin{eqnarray*}
 \partial^{\alpha} \tau^j =
\sum_{k = 1,...,m; \vert \beta \vert \leq q+r-1}
(b^{j}_{k\beta}(y) \partial^{\beta} \tau^k), j = 1,...,m, \vert \alpha \vert = q +r
\end{eqnarray*}

This implies by the chain rule and reccurence that all derivatives of
$\tau^j$
of order $\geq q +r$ at $y^0$ are determined by derivatives of order
$\leq q+r-1$, which means that the dimension of $Sol({\cal R}_q)$ is finite.

This proof is quite constructive  and allows to obtain explicit recurcive
formulae for the Taylor expansions at $y^0$ of solutions of $({\cal R}_0)$.
This also means that the dimension of $Sol({\cal R}_q)$ is majorated by
$dim J^t(n,m)$ where $d = type({\cal R}_q) - 1$. Of course this estimate is
not
precise since the partial  derivatives at $y^0$ of $\tau$ of order $\leq d$ satisfy
a system of linear algebraic equations $({\cal L})$ arising from the 
equations of $({\cal R}_{q+r})$ of order $< (q+r)$. Solving
this
system we can presisely determine the dimension of the space $Sol({\cal
R}_q)$
for any concrete system $({\cal R}_q)$. More precisely, applying
 the Cramer rule to $({\cal L})$ we 
can represent some partial derivatives of $\tau$ at $y^0$ of order
$\leq d$ ({\it principal derivatives}) as linear combinations of
others ({\it parametric derivatives}). The number of parametric
derivatives is equal to the dimension of $Sol({\cal R}_q)$  and they
form a set of natural parameters on $Sol({\cal R}_q)$.

Let  $({\cal R}_q^{\varepsilon})$
be an analytic family of linear systems given by

\begin{eqnarray*}
({\cal R}_q^{\varepsilon}): \sum_{j = 1,...,m; \vert \alpha \vert \leq q}
a^{i}_{j\alpha}(\varepsilon,y) \partial^{\alpha} \tau^j = 0, i= 1,...,s
\end{eqnarray*}
where $a^{i}_{j\alpha}$ are holomorphic functions in  $y$ and 
real analytic in $\varepsilon$, 
with $\varepsilon$ being in a neighborhood of the origin in $\R^k$. The
following obvious observation turns out to
be very useful:

\begin{e-pro}
Suppose that the system ${\cal R}_q^0$ is of finite type. Then for every
$\varepsilon$ close enough to the origin the system $({\cal
R}_q^{\varepsilon})$
is of finite type and $type({\cal R}_q^{\varepsilon}) \leq type({\cal
R}_q^0)$.
Furthermore, $dim Sol({\cal R}_q^{\varepsilon}) \leq dim Sol({\cal R}_q^0)$.
\end{e-pro}

The proof is immediate since the rank of a linear algebraic system 
defining the symbol of the prolonged system does not decrease with respect to small perturbations of the coefficients so  $type({\cal R}_q^{\varepsilon})
 \leq type({\cal R}_q^0)$. Similarly, if $({\cal L}^{\varepsilon})$ is a 
linear algebraic system for the partial derivatives of order $<
type({\cal R}_q)$ arising from the equations of the lower orders, then
$rank ({\cal L}^{\varepsilon}) \geq rank ({\cal L}^0)$ and the 
number of the parametric derivatives decreases so  
$dim Sol({\cal R}_q^{\varepsilon}) \leq dim Sol({\cal R}_q^0)$.

In general a linear system of order $q$ may contain some equations of order
$< q$.
However, if we add to such a system all the equations of order $\leq q$
obtained from
the equations of lower order by taking all the partial derivatives of
a  suitable order, we obtain a system with the same space of solutions. We
call such a system the {\it completion} of $({\cal R}_q)$ or  
the {\it completed} system  $({\cal R}_q)$. We also point out that
every linear system can be reduced to a system of the first order by
introducing the supplementary dependent variables; so one may work
with these systems only.

Applying these results to the completed Lie equations 
deduced in the previous section for an involutive 
system $({\cal S}^0)$ and its holomorphic involutive deformation,
 we obtain the following

\begin{e-theo}
Suppose that the completed  Lie equations for $({\cal S}^0)$ form a system of finite type $d$ at some point
$(x^0,u^0)$.  Then $dim Lie({\cal S}^0)$ is finite and for any $\varepsilon$
close enough to the origin $dim Lie({\cal S}^{\varepsilon}) \leq dim Lie(
{\cal S}^0)$.
\end{e-theo}

In view of this result it is of clear interest the question how to check up
if a given system is of finite type. On of the possibilities here is to
consider its {\it characteristic variety}. Let $\lambda$ be a vector of
$\cc^n$. We use the notation $\lambda^{\alpha} = \lambda^{\alpha_1}...
\lambda^{\alpha_n}$. A vector $\lambda$ is called a {\it characteristic}
(co)vector at $y$ if the linear map $\sigma_{\lambda}(y): \cc^m \longrightarrow \cc^s$
given by the matrix $\sigma_{\lambda}(y): \sum_{ \vert \alpha \vert = q}
a^{i}_{j\alpha}(y) \lambda^{\alpha}$
is not injective. The set of of such $\lambda$ is an algebraic variety in
$\cc^n$ which is called the {\it characteristic variety} at $y$ and is
denoted by $Char_y({\cal R}_q)$.

The following criterion is useful (see \cite{Po}, p.195): a system $({\cal R}_q)$ is of finite type if and only if $Char_y({\cal
R}_q)$
is zero for every $y$ (we do not use it in the present paper).

Of course, this statement says nothing about a value of the type of $({\cal
R}_q)$.
However, if the system $({\cal R}_q)$ is known to be of finite type, its
type can be determined by direct computations using the study of a finite
number of prolongations and their symbols, i.e. by means of the elementary
linear algebra tools.

As an example we study the Lie equations in the simplest classical case of
a second order ordinary differential equation.

We denote by $x \in \cc$ and $u \in \cc$ the independent and dependent
variables
respectively and consider a holomorphic equation $({\cal S}): u_{xx} = F(x,u,u_x)$. This equation define a hypesurface in the jet space $J^2(1,1)$ :
$({\cal S}_2): u_{11} = F(x,u,u_1)$.

A holomorphic vector field $X = \theta \frac{\partial}{\partial x} + \eta \frac{\partial}{\partial u}$
is an infinitesimal symmetry of $({\cal S})$ if and only if
its 2-prolongation $ X^{(2)} = X + \eta_1\frac{\partial}{\partial u_1} +
\eta_{11}\frac{\partial}{\partial u_{11}}$
is tangent to $({\cal S}_2)$ that is $X^{(2)}(u_{11} - F(x,u,u_1)) = 0, (x,u,u_1,u_{11}) \in ({\cal S}_2)$.

The coefficients have the following expessions:

\begin{eqnarray*}
& &\eta_1 = \eta_{ x} + \left (
\eta_{ u} - \theta_{ x} \right ) u_1 -
\theta_{ u} (u_1)^2,\\
& &\eta_{11} = \eta_{ xx} +
\left ( 2\eta_{ x  u} -
\theta_{ xx} \right ) u_1 +
\left (  \eta_{ uu} - 2 \theta_{ x  u} \right ) (u_1)^2 - 
 \theta_{ uu}(u_1)^3 +
\left (  \eta_{ u} - 2 \theta_{ x} \right ) u_{11} - 
3\theta_{ u}u_1 u_{11}
\end{eqnarray*}

Consider the expansion $F(x,u,u_1) = \sum_{\nu \geq 0} f^{\nu}(x,u)
(u_1)^{\nu}$;
after elementary computations following the decribed above general method
we obtain the following system $({\cal R}_2)$ of infinitesimal Lie
equations:

\begin{eqnarray*}
& &\eta_{xx} = 2 f^0 \theta_{x} + f^1\eta_x - f^0\eta_u + f^0_x\theta +
f^0_u\eta,\\
& &2\eta_{xu}  - \theta_{xx} = f^1\theta_x - 3f^0\theta_{u} + f^1_x\theta +
f^1_u\eta,\\
& &\eta_{uu} - 2\theta_{xu} = 2f^1\theta_u + 3f^3\eta_x + f^2_x\theta +
f^2_u\eta,\\
& &-\theta_{uu} = -f^3\theta_x + f^2\theta_u + 4f^4\eta_x + f^3_x\theta +
f^3_u\eta,\\
& &(2 - \nu)f^{\nu}\theta_x + (4 - \nu)f^{\nu - 1} \theta_u + (\nu +
1)f^{\nu + 1}
\eta_x + f^{\nu}_x\theta + f^{\nu}_u\eta = 0, \nu \geq 4
\end{eqnarray*}

Actually only a finite number of these equations are independent. But we show
that the first 4 second order equations form a finite type system. Thus,
we consider a system $({\cal R}_2^{'})$ :

\begin{eqnarray}
\label{55.1}
\eta_{xx} = 2 f^0 \theta_{x} + f^1\eta_x - f^0\eta_u + f^0_x\theta +
f^0_u\eta
\end{eqnarray}
\begin{eqnarray}
\label{55.2}
2\eta_{xu}  - \theta_{xx} = f^1\theta_x - 3f^0\theta_{u} + f^1_x\theta +
f^1_u\eta
\end{eqnarray}
\begin{eqnarray}
\label{55.3}
\eta_{uu} - 2\theta_{xu} = 2f^1\theta_u + 3f^3\eta_x + f^2_x\theta +
f^2_u\eta
\end{eqnarray}
\begin{eqnarray}
\label{55.4}
-\theta_{uu} = -f^3\theta_x + f^2\theta_u + 4f^4\eta_x + f^3_x\theta +
f^3_u\eta
\end{eqnarray}

The symbol $G'_2$ of this system is a linear 2- dimensional subspace of the
space $\cc^6$ with coordinates $v^1_{11}$, $v^1_{12}$, $v^1_{22}$,
$v^2_{11}$, $v^2_{12}$, $v^2_{22}$ defined by the equations

\begin{eqnarray*}
v^2_{11} = 0, 2v^2_{12} - v^1_{11} = 0, v^2_{22} - 2v^1_{12} = 0, v^1_{22} =
0
\end{eqnarray*}

A vector $\lambda \in \cc^2$ will be characteristic if and only if the
matrix with the lines $(0 ,\lambda_1^2)$, $(\lambda_1^2
,2\lambda_1\lambda_2)$,
$(-2\lambda_1\lambda_2,\lambda_2^2)$, $(\lambda_2^2,0)$ has the rank $\leq
1$; this implies the the characteristic variety is equal
to zero and so our system is of finite type.

Its 1-prolongation $G'_3$ is a subspace of $\cc^8$ with the coordinates
$v^1_{111}$, $v^1_{112}$, $v^1_{122}$, $v^1_{222}$,
$v^2_{111}$, $v^2_{112}$, $v^2_{122}$, $v^2_{222}$ given by the equations

\begin{eqnarray*}
& &v^2_{111} = 0, v^2_{112} = 0, v^1_{122} = 0, v^1_{222} = 0,\\
& &2v^2_{112} - v^1_{111} = 0, 2v^2_{122} - v^1_{112} = 0,\\
& &v^2_{122} - 2v^1_{112} = 0, v^2_{222} - 2v^1_{122} = 0
\end{eqnarray*}
so we see immediately that $G'_3 = \{ 0 \}$, i.e. $({\cal R}'_2)$ is
of type 1. Solving its 1-prolongation $({\cal R}'_3)$ with respect to
the partial derivatives of the third order, we obtain the following
explicit representations:

\begin{eqnarray*}
& &\theta_{xxx} = -f^1\theta_{xx} + 7f_0\theta_{xu} + 2f^1\eta_{xu} -
2f^0\eta_{uu} + 4(f^0_u - f^1_x - f^1_x)\theta_x + 5f^0_x\theta_u +
f^1_u\eta_x \\
& &+ (2f^0_{xu} - f^1_{xx})\theta +
(2f^0_{uu} - f^1_{xu})\eta,\\
& &\theta_{xxu} = -f^0\theta_{uu} - f^1\theta_{xu} - 2f^3\eta_{xx} +
(1/3)(f^1_u - 2f^2_x)\theta_x - (f^0_u + f^1_x)\theta_u -
(1/3)(5f^3_x + 2f^2_u)\eta_x\\
& &+ (1/3)f^1_u\eta_u + (1/3)(f^1_{xu} - 2f^2_{xx})\theta +
(1/6)(2f^1_{uu} - f^2_{xu})\eta,\\
& &\theta_{xuu} = f^3\theta_{xx} - 4f^4\eta_{xx} - f^2_x\theta_u - (4f^4_x +
f^3_u)\eta_x - f^3_{xx}\theta - f^3_{xu}\eta,\\
& &\theta_{uuu} = f^3\theta_{xu} - f^2\theta_{uu} - 4f^4\eta_{xu}
+ f^3_u\theta_x - (f^2_u + f^3_x)\theta_u - 4f^4_u\eta_x
- f^3_u\eta_u - f^3_{xu}\theta - f^3_u\eta,\\
& &\eta_{xxx} = 3f^0_x\theta_x + (f^1_x + f^0_u)\eta_x - f^0_x\eta_u
+ 2f^0\theta_{xx} + f^1\eta_{xx} - f^0\eta_{xu} +
f^{0}_{xx}\theta + f^0_{xu}\eta,\\
& &\eta_{xxu} = 2f^0\theta_{xu} + f^1\eta_{xu} - f^0\eta_{uu} +
2f^0_u\theta_x
+ f^0_x\theta_u + f^1_u\eta_x + f^0_{xu}\theta + f^0_{uu}\eta,\\
& &\eta_{xuu} = - 2f^0\theta_{uu} - f^3\eta_{xx} +(1/3)(2f^1_u -
 f^2_x)\theta_x
 -2f^0_u\theta_u - (1/3)(f^3_x + f^2_u)\eta_x + (2/3)f^1_u\eta_u\\
& &+ (1/3)(2f^1_{xu} - f^2_{xx})\theta + (1/3)(2f^1_{uu} - f^2_{xu})\eta,\\
& &\eta_{uuu} = 2f^3\theta_{xx} + 2f^1\theta_{uu} - 8f^4\eta_{xx} +
 3f^3\eta_{xu}
 - (2f^1_u - f^2_x)\theta_u + (f^3_u - 8f^4_x)\eta_x + f^2_u\eta_u\\
& &+ ( f^2_{xu} - 2f^3_{xx})\theta + (f^2_{uu} - 2f^2_{xu})\eta
\end{eqnarray*}

Fix a point $(x_0,u_0)$ and attach the values $a_1: = \theta(x_0,u_0)$,
$a_2: = \eta(x_0,u_0)$, $a_3:= \theta_x(x_0,u_0)$, $a_4:=
\theta_u(x_0,u_0)$, $a_5:= \eta_x(x_0,u_0)$, $a_6:= \eta_u(x_0,u_0)$,
$a_7 = \theta_{xx}(x_0,u_0)$, $a_8 = \theta_{xu}(x_0,u_0)$ to the 
parametric derivatives. Then the values
of all second order derivatives of $\theta$, $\eta$ at $(x_0,u_0)$ are
determined by (\ref{55.1})- (\ref{55.4}) and the values of all derivatives at $(x_0,u_0)$
of order $\geq 3$ are determined by the former expresions for the third
order partial derivatives via the chain rule. This means
that $dim Lie({\cal S}) \leq 8$ and this estimate is precise since
in the flat case where $F \equiv 0$ one has $dim Lie({\cal S}) = 8$.

Of course, the constructed vector fields are in general just the {\it
candidates} to
be in $Lie({\cal S})$ since we still have additional first order equations
in
the Lie equations $({\cal R}_q)$. The fact that $\theta$, $\eta$ satisfy
these equations imposes additional analytic restrictions on the parameters
$a_j$ so actually $Lie({\cal S})$ is parametrized by a some analytic
subvariety
in the space $\cc^8$ of the parameters $a_j$.

The present description of symmetries of a second order ordinary 
differential equation has been obtained by L.Dickson \cite{Di}.
Since the Segre family of  a Levi nondegenerate hypersurface in 
$\cc^2$ is a set of solutions of such equation, the present method 
allows to obtain an explicit parametrization of its automorphism group.
 This argument
can be directly generalized to second order holomorphic involutive PDE symmetries

\begin{eqnarray*}
u^k_{x_ix_j} = F^k_{ij}(x,u,u_x), k = 1,...,m, i,j =1,...,n
\end{eqnarray*}
Using this method and the explicit formulae for the 2-prolongation of a 
vector field on $\cc^n \times \cc^m$, the author proved in \cite{Su}
that the Lie algebra of infinitesimal symmetries of 
such a system has a dimension $\leq (n+m+2)(n+m)$ and every infinitesimal 
symmetry is determined by a second order Taylor expansion at a given point 
(the Lie equations are of type 1). In the special case where $n = 1$ i.e. 
for a system of ordinary differential equations this result was established 
 by F.Gonzales-Gascon and A.Gonzales-Lopez \cite{GL} (see also \cite{Ol2}).  In particular, this implies the results of Tanaka \cite {Ta} and
Chern - Moser \cite{CM} on the majoration of the dimension of the automorphism 
group of a real analytic Levi nondegenerate hypersurface in $\cc^{n+1}$, 
its parametrization etc.

It is important to emphasize that such an explicit parametrization of the
Lie algebra of infinitesimal symmetries can be obtained for every system
with the Lie equations of finite type. In what follows we restrict ourselves
just by the study of symbols of the Lie equations in order to avoid
complicated formulae.

We conclude this section by a  statement concerning the
special case of linear PDE systems with constant coefficients. The main
example
of these systems is given by the Lie equations for a flat manifold derived
in
the previous section.

Consider a linear PDE system with {\it constant} coefficients of the form

\begin{eqnarray*}
({\cal R}_q): \sum_{i,\vert \alpha \vert =  q_k} a^k_{i\alpha}
\partial^{\alpha}
u^i = 0, k = 1,...,K
\end{eqnarray*}
where $q_k =max_k q_k$.
{\it We emphasize that every equation of this system  of order
$q_k$ contains the partial derivatives of the same order $q_k$
only}. In particular, the Lie equations for a flat system deduced in
the 
previous section are of this class.

A holomorphic in a neighborhood of the origin map $u = (u^1,...,u^m)$ is a
solution of $({\cal R}_q)$ if and only if

\begin{eqnarray*}
\partial^{\beta}(\sum a^k_{i\alpha}\partial^{\alpha}u^i)\vert_{x=0} =
\sum a^k_{i\alpha}\partial^{\beta + \alpha}u^i\vert_{x = 0} = 0,
k = 1,...,K
\end{eqnarray*}
for every $\beta$.

This is equivalent to

\begin{eqnarray}
\label{5.1}
\sum_{i;\vert \alpha \vert = q_k, \vert \beta \vert = s- q_k}
a^k_{i\alpha} (\partial^{\beta + \alpha}u^i\vert_{x=0}) = 0, k = 1,...,K,
s = q, q+1,...
\end{eqnarray}

In the complex affine space with the coordinates $(v^i_{i_1...i_s})$, $i\in
\{ 1,..., m \}$, $i_1 \leq ... \leq i_s$, $i_j \in \{ 1,..., n \}$ consider
a subspace $V_s$ defined by the linear algebraic system

\begin{eqnarray*}
\sum_{i;\vert \alpha \vert = q_k, \vert \beta \vert = s- q_k}
a^k_{i\alpha} v^i_{\beta + \alpha} = 0, k = 1,...,K
\end{eqnarray*}
for $s = q, q+1,...$.

\begin{e-pro}
\label{pro5.5}
The dimension of the space $Sol({\cal R}_q)$ is finite if and only if
there exists an $s$ such that $V_s = \{ 0 \}$.  In this case the 
completion of $({\cal R}_q)$ is a 
system is of finite type and every
solution is a polynomial of degree $< s$.
\end{e-pro}

\proof Suppose that there exists an $s$ such that $V_s = \{ 0 \}$. In view
of
(\ref{5.1}) this means that the completion of $({\cal R}_q)$  is a system
of finite type majorated by $s$. Moreover, (\ref{5.1}) shows that
in this case all partial derivatives of $u$ of order $s$ vanish identically.

Let now the dimension of $Sol({\cal R}_q)$ is finite. Suppose by
contradiction
that there exists an increasing sequence $(s_t)$ such that every $V_{s_t}$
is
non-trivial. Let  $(v^i_{i_1...i_{s_t}})$ be a non-zero vector in $V_{s_t}$.
Consider the map $u_t = (u^1_t,...,u^m_t)$ whose components are the
homogeneous polynomials of degree $s_t$ satisfying
$\frac{\partial^{s_t}u^i_t}{\partial x_{i_1}...\partial
x_{i_{s_t}}}(0)  = v^i_{i_1...i_{s_t}}$. Then for
every $t$ the function $u_t$ satisfies  (\ref{5.1}) for $s = s_t$; but since it
is homogeneous polynomial of degree $s_t$, clearly it satisfies (\ref{5.1}) for all other
$s$.
Therefore, every $u_t$ is a solution of $({\cal R}_q)$: a contradiction.

 In particular, we have the following

 \begin{e-cor}
\label{cor5.6}
 Suppose that $({\cal R}_q)$ has a finite dimensional solution space and let
 $({\cal R}_q^{\varepsilon})$ be its  holomorphic deformation.
Then
 for every $\varepsilon$ small enough $dim Sol({\cal R}_q^{\varepsilon})
 \leq dim Sol({\cal R}_q)$.
 \end{e-cor}

\section{ General flat systems with first order linear relations}

In this section we consider a flat system $({\cal S})$ of the form

\begin{eqnarray*}
& & u^1_{x_ix_j} = 0, i,j = 1,...,n \\
& & u^k_x = A^k u^1_x, k = 2,...,n
\end{eqnarray*}
with $n$ independent and $m$ dependent variables. We apply a geometric 
method in order to describe the symmetries of this system without 
computations. The basic idea goes back to S.Lie - G.Scheffers \cite{LS}
(see also \cite{GL}); a related result also was obtained by B.Shiffman
\cite{Sh}. The present proof is a direct generalization of  author's 
argument about the rationality of holomorphic maps between quadrics 
in $\cc^n$ \cite{Su2}.

\begin{e-theo}
Suppose that the matrices $A^1: = Id_n$, $A^2$,...,$A^m$ are linearly
independent. Then $Lie({\cal S})$ is finite dimensional.
\end{e-theo}

\proof Fix an infinitesimal symmetry $X \in Lie({\cal S})$ and for
$t \in \cc$ close enough to the origin  consider the flow $f(t,x,u) = e^{tX}$
generated by $X$.

The set $Sol({\cal S})$ of solutions of $({\cal S})$ is an $(n+m)$- parameter family of affine
subspaces of $\cc^n \times \cc^m$ of the form 
$Q(\zeta,\omega) = \{ (x,u): u = \omega + <x,A\zeta> \}$
where $\omega + <x,A\zeta>  = \omega^j + <x,A^j>$ , $j = 1,...,m$. The
parameters $(\zeta,\omega) \in \cc^{n+m}$ give a natural
holomorphic coordinate
 system on $Sol({\cal S})$ which is an $(n+m)$-dimensional complex 
manifold.

The fact that $f_t$ takes any solution to another solution means that
for any $(\zeta,\omega)$ there exists a point $(\zeta^*_t, \omega^*_t)$
such that $f_t(Q(\zeta,\omega)) = Q(\zeta^*_t,\omega^*_t)$ that is

\begin{eqnarray}
\label{6.1}
h_{t}(x,\omega + <x,A\zeta>) = \omega_t^* +
<g_t(x,\omega + <x,A\zeta>), A\zeta^*_t>
\end{eqnarray}
where $f_t = (g_t,h_t)$.

Thus, $f_t$ induces a map

\begin{eqnarray*}
& & f^*_t: Sol({\cal S}) \longrightarrow Sol({\cal S}),\\
& & f^*_t: (\zeta,\omega) \mapsto (\zeta^*_t, \omega^*_t)
\end{eqnarray*}

\begin{e-lemme}
The family $\{ f_t \}$ is a family of biholomorphisms holomorphically
depending on the parameter $t$.
\end{e-lemme}

\proof The image $f_t(Q(\zeta,\omega))$ is given by

\begin{eqnarray*}
\{ (x^*,u^*):  (x^*,u^*) = (g(t,x,\omega + <x,A\zeta>),
h(t,x,\omega + <x,A\zeta>)), x \in \cc^n \}.
\end{eqnarray*}

For $t = 0$ one has
$(g_0(\bullet),h_0(\bullet)) = (x,u)$ so for $t$ small enough the
implicit function theorem can be applied to $x^* = g(t,x,\omega +
<x,A\zeta>$ and $x = x(t,x^*,\zeta,\omega)$ is holomorphic. Substituting it
to $u^* = h(t,x,\omega + <x,A\zeta>)$ we obtain
$u^* = \varphi(t,x^*,\zeta,\omega)$
and $\varphi$ is holomorphic. On the other hand,
$f_t(Q(\zeta,\omega)) = Q(\zeta^*_t,\omega^*_t)$ so 
$\varphi(t,x^*,\zeta,\omega) = \omega^*_t + <x^*,A\zeta^*_t>$. In particular,
$\varphi_1(t,x^*,\zeta,\omega) = \omega_1^* + x_1^*(\zeta_t^*)_1 + ...+
x_n^*(\zeta^*_t)_n$
so every $\zeta_j^* = \zeta_j^*(t,\zeta,\omega)$ is holomorphic and
obviously
$\omega^* = \omega^*(t,\zeta,\omega)$ is holomorphic.

Consider the vector fileds ${\cal L}_{\nu} = \frac{\partial}{\partial \zeta_{\nu}} -
\sum_{k = 1}^{m}\left ( \sum_{j=1}^n a^k_{i\nu} x_i \right )
\frac{\partial}{\partial \omega^k}$
where $A^k = (a^k_{ij})$.

Applying them to (\ref{6.1}) we get

\begin{eqnarray}
\label{6.2}
{\cal L}_{\nu}((\omega^*_{j})_t) + <g_t(x,\omega + <x,A\zeta>),{\cal
L}_{\nu}A^j\zeta_t^*> = 0
\end{eqnarray}
Consider $(\ref{6.1})$, $(\ref{6.2})$ as a linear system with respect to
components of $f_t$. Since $(\zeta_0^*,\omega_0^*) \equiv (\zeta,\omega)$,
this
system contains an $(n +m) \times (n+m)$ subsystem with the determinant $\neq 0$ for
$t$ small enough. Applying the Cramer rule we obtain that for any
$(t,\zeta,\omega)$
fixed the map $f_t(x,\omega + <x,A\zeta>)$ is a rational map in $x$. Moreover,
the degree of every such a map is uniformly bounded by $n$.

The last step of the proof is to show the the space of solutions $({\cal
S})$
is "large enough".

Set $(e_k = (0,...,1,...,0) \in \cc^n$ ($1$ on the $k$-position) and
consider
the vectors $v_k(\zeta) = (e_k,<e_k,A^1\zeta>,...,<e_k,A^m\zeta>)$
(so $v_k(\zeta) \in Q(\zeta,0)$).

\begin{e-lemme}
The linear hull of $\{ v_k(\zeta), \zeta \in \cc^n \}$ coincides with
$\cc^n$.
\end{e-lemme}
\proof If the statement is false, there exists a $\lambda \in
\cc^{n+m} \backslash \{ 0 \}$ such that $<\lambda,v_k(\zeta)> = 0$ for
any $k,\zeta$ that is
$\lambda_k + \lambda_{n+}<e_k,A^1\zeta> + ...+ \lambda_{n+m}<e_k,A^m\zeta> =
0$ for all $\zeta$, $k$; therefore $\lambda_k = 0$ for every $k= 1,...,n$ and
so $<e_k,(\lambda_{n+1}A^1 + ...+ \lambda_{n+m}A^m)\zeta> = 0$
for every $k$, $\zeta$, that is $\lambda_{n+1}A^1+...+\lambda_{n+m}A^m = 0$
:
a contradiction which proves the lemma.

Fix now $(n+m)$ linearly independent complex lines $l^1,...,l^{n+m}$, every
$l^j$ is in some $Q(\zeta^j,0)$ through the origin. Every line generates a
family of parallel lines and any line of such a family is in
$Q(\zeta^j,\omega)$
for some $\omega$. After a linear change of variables in $\cc^{n+m}$ these
families become the coordinate ones and the classical separate rationality
theorem \cite{BM} implies that $f_t$ is a rational map of degree $\leq n$ for any
$t$ small enough that is

\begin{eqnarray*}
f_t(x,u) = \frac{\sum_{\vert I \vert = 0}^n a_I(t)(x,u)^I}
{\sum_{\vert I \vert = 0}^n b_J(t) (x,u)^J}.
\end{eqnarray*}
Hence, $X = \frac{df_t}{dt} \vert \{ t = 0 \}$ is a vector field with
rational
coefficients of degree $\leq n^2$. Every such a coefficient is uniquely
determined by a finite number $d = d(n^2)$ of terms of its Taylor expansion
at the origin. Therefore, the dimension of $({\cal S})$ is finite.
This completes the proof of the theorem.

We say that a flat system $({\cal S})$ is {\it nondegenerate} if it
satisfies the hypothesis of our proposition that is the matrices $A^j$ are
linearly independent.

From Proposition \ref{pro5.5} we obtain the following 

\begin{e-cor}
The completed Lie equations of a nondegenerate flat  system
 $({\cal S})$ form a PDE system of finite type  and 
 every infinitesimal symmetry $X \in Lie({\cal S})$ has polynomial 
coefficients of uniformly bounded degree.
\end{e-cor}

Corollary \ref{cor5.6} implies now one of our  main results:

\begin{e-theo}
If $({\cal S}^{\varepsilon})$ is an involutive holomorphic deformation of
a nondegenerate flat system $({\cal S})$,
then $dim Lie({\cal S}^{\varepsilon}) \leq dim Lie({\cal S})$.
\end{e-theo}

Now we can apply the obtained results in order to study biholomorphisms of
Cauchy-Riemann manifolds.

Let ${\cal M}$ be a generic real analytic Levi nondegenerate submanifold
in $\cc^{n+m}$ through the origin. After a biholomorphic change 
of coordinates it can be represented  in the form $w + \overline{w} =
<L(z),\overline{z}> + o(\vert Z \vert^2)$. Denote by ${\cal M}_{flat}$
the corresponding quadric: $w + \overline{w} =
<L(z),\overline{z}>$. For real $\varepsilon$ 
close enough to the origin  consider the 
following change of variables: $z = \varepsilon  z', w = \varepsilon^2  w'$.

In the new coordinates (we omit the primes) we get the manifold ${\cal M}^{\varepsilon}: w + \overline{w} = 
<L(z),\overline{z}> + (1/\varepsilon^2)R(\varepsilon z,
\varepsilon \overline z, \varepsilon^2 w, \varepsilon^2 \overline w)$
biholomorphic to ${\cal M}$ for every $\varepsilon$. Since the 
function $(1/\varepsilon^2)R(\varepsilon z,
\varepsilon \overline z, \varepsilon^2 w, \varepsilon^2 \overline w)$
extends to a function real analytic in $\varepsilon$ in a neighborhood
of the 
origin and vanishing at the origin, the 
system ${\cal S}({\cal M}^{\varepsilon})$ defining the Segre family of
${\cal M}^{\varepsilon}$ is 
a holomorphic involutive deformation of the flat system defining the Segre 
family of ${\cal M}_{flat}$.

It follows from the results of the previous sections
that we have established the following  result:

\begin{e-cor}
$Aut({\cal M})$ is a finite dimensional real Lie group. 
Moreover, $dim Aut({\cal M})$ is majorated by the
complex dimension of the flat PDE system defining 
the Segre family of ${\cal M}_{flat}$. 
\end{e-cor}

Various results of this type for this and  more general classes of CR
manifolds have been obtained by several authors 
\cite{BER2,BER3, Be1,Be2, EIS,Lo,St,Tu,Za}
using  different methods. 
We emphasize that our method can be adapted to a much more general situation 
and allows to obtain  many additional information on the structure of the automorphism group.

Remark. We have introduced the small parameter $\varepsilon$ by
analogy with  the well-known scaling techniques (see for instance
\cite{BP2}). On the other hand, in our situation this argument can be 
considered as an application of the general PDE method of small
parameter widely known in the classical mechanics. 

The geometric method employed in this section   allows to obtain only  
an inprecise estimate of the type of the Lie equations. In order 
to determine
this type precisely, a direct linear algebra computations can be used. In the next section we consider the 
special case of system with two dependent and two independent variables
and show how the computations of the type can effectively be done.

\section{ Flat systems with linear relations, case $n = 2$, $m = 2$}

In the present section we consider the special case  of study of
(infinitesimal)
symmetries of flat systems with first order relations.

   Consider the following flat system $({\cal S})$
given by

\begin{eqnarray*}
& & u^j_{x_1x_1} = 0, u^j_{x_1x_2} = 0, u^j_{x_2x_2} = 0, j = 1,2\\
& & u^2_{x_1} = a_{11}u^{1}_{x_1} + a_{12}u^{1}_{x_2},\\
& & u^2_{x_2} = a_{21}u^{1}_{x_1} + a_{22}u^1_{x_2}\\
\end{eqnarray*}

Our goal is to establish the following

\begin{e-pro}
Suppose that the matrcies $Id_2$, $A$ are linearly independent that is
$({\cal S})$ is nondegenerate. Then
 the corresponding Lie equations of $({\cal S})$ form a PDE system
of finite type 1.
\end{e-pro}

Let a holomorphic vector field $X = \theta^1 \frac{\partial}{\partial x_1} + \theta^2
\frac{\partial}{\partial x_2}
+ \eta^ 1 \frac{\partial}{\partial u^1} + \eta^2 \frac{\partial}{\partial
u^2}$
be in $Lie({\cal S})$. First and second prolongations are

\begin{eqnarray*}
X^{(1)} = X + \eta^1_1 \frac{\partial}{\partial u^1_1} +
\eta^1_2 \frac{\partial}{\partial u^1_2} +
\eta^2_1 \frac{\partial}{\partial u^2_1}
+ \eta^2_2 \frac{\partial}{\partial u^2_2}
\end{eqnarray*}

\begin{eqnarray*}
& & X^{(2)} = X^{(1)} + \eta^1_{11} \frac{\partial }{\partial u^1_{11}} +
\eta^1_{12} \frac{\partial }{\partial u^1_{12}}
+ \eta^1_{22} \frac{\partial }{\partial u^1_{22}}+ \eta^2_{11}
\frac{\partial }{\partial u^2_{11}} +
\eta^2_{12} \frac{\partial }{\partial u^2_{12}}
+ \eta^2_{22} \frac{\partial }{\partial u^2_{22}}
\end{eqnarray*}

Following the general method described above, we have to consider the first
order Lie equations:

\begin{eqnarray*}
& &\eta^2_1 \vert ({\cal S})^{(2)} = a_{11}\eta^1_1 \vert ({\cal S})^{(2)} +
a_{12}\eta^1_2 \vert ({\cal S})^{(2)} \\
& &\eta^2_2 \vert ({\cal S})^{(2)} = a_{21}\eta^1_1 \vert ({\cal S})^{(2)} +
a_{22}\eta^1_2 \vert ({\cal S})^{(2)}
\end{eqnarray*}

Computing the restrictions $\eta^2_1 \vert ({\cal S})^{(2)}$ and
comparing the coefficients near the powers of $u^k_j$, we obtained the
following linear first order PDE systems  with constant coefficients for
$\theta$, $\eta$:

\begin{eqnarray*}
  \eta^2_{x_1} = a_{11}\eta^1_{x_1} + a_{12}\eta^1_{x_2},
\eta^2_{x_2} = a_{21}\eta^1_{x_1} + a_{22}\eta^1_{x_2}
\end{eqnarray*}
and
\begin{eqnarray*}
& & \eta^2_{u^1} + a_{11}\eta^2_{u^2} =
a_{11}\eta^1_{u^1} + (a^2_{11} + a_{12}a_{21})
\eta^1_{u^2} - a_{12}\theta^1_{x_2},\\
& & a_{12}\eta^2_{u^2} =
a_{12}\eta^1_{u^1} + a_{12}(a_{11} + a_{22})
\eta^1_{u^2} +
a_{12}\theta^1_{x_1} + (a_{11} + a_{22})
\theta^2_{x_1}
- a_{12}\theta^2_{x_2},\\
& & a_{21} \eta^2_{u^2} =
a_{21}\eta^1_{u^1} + (a_{21}a_{11} + a_{22}a_{21})
\eta^1_{u^2} -
a_{11}\eta^1_{x_2} -
a_{21}\theta^1_{x_1} -
a_{22}\theta^1_{x_2},\\
& &\eta^2_{u^1} +
a_{22}\eta^2_{u^2} =
(a_{21}a_{12} + a^2_{22})\eta^1_{u^2} +
a_{22}\eta^1_{u^1} +
a_{12}\theta^1_{x_2} -
a_{21}\theta^2_{x_1}
\end{eqnarray*}

In view of our condition of linear independence of $Id$, $A$ this last
system implies that

\begin{eqnarray*}
\eta^2_{ u^1} =
\phi_1( \eta^1_{ u^j}, \theta^i_
{ x_k}), \eta^2_{ u^2} =
\phi_1( \eta^1_{ u^j}, \theta^i_
{ x_k})
\end{eqnarray*}
where $\phi_s$ are linear functions.

Finally, we have two series of equations:

\begin{eqnarray*}
& & a_{21}\left (\theta^2_{u^1} +
a_{11} \theta^2_{u^2} -
a_{12}\theta^1_{u^2}  \right ) = 0, \\
& & (a_{11} - a_{22})  \left ( \theta^2_{u^1}-
a_{12}\theta^1_{u^2} +
a_{11}\theta^2_{u^2} \right ) = 0, \\
& & a_{12} \left ( \theta^2_{u^1} -
a_{12}\theta^1_{u^2} +
a_{11}\theta^2_{u^2} \right ) = 0
\end{eqnarray*}

\begin{eqnarray*}
& &a_{21}\left ( \theta^1_{u^1} -
a_{21} \theta^2_{u^2} -
a_{22}\theta^1_{u^2}  \right ) = 0 \\
& &(a_{11} - a_{22}) \left ( \theta^1_{u^1} +
a_{22}\theta^1_{u^2} -
a_{21}\theta^2_{u^2} \right ) = 0\\
& &a_{12} \left ( \theta^1_{u^1} +
a_{22}\theta^1_{u^2} -
a_{21}\theta^2_{u^2} \right ) = 0
\end{eqnarray*}

In view of the linear independence of the matrices $Id_2$, $A$ this implies
that

\begin{eqnarray*}
\theta^1_{u^1} = -a_{22}\theta^1_{u^2} +
a_{21}\theta^2_{u^2},
\theta^2_{u^1} = a_{12}\theta^1_{u^2} -
a_{11}\theta^2_{u^2}
\end{eqnarray*}

 It is useful to consider the  differential consequences of these equalities:

 \begin{eqnarray*}
& & \theta^1_{u^1u^2} = -a_{22}\theta_{u^2u^2} + a_{21}\theta_{u^2u^2},\\
& &\theta^2_{u^1u^2} = a_{12}\theta^1_{u^2u^2} - a_{11}\theta_{u^2u^2},\\
& &\theta^1_{u^1u^1} = (a_{12}a_{21} +a^2_{22})\theta^1_{u^2u^2} -
(a_{22}a_{21} + a_{21}a_{11})\theta^2_{u^2u^2},\\
& &\theta^2_{u^1u^1} = -(a_{12}a_{22} + a_{11}a_{12})\theta^1_{u^2u^2} +
(a_{12}a_{21} + a^2_{11})\theta^2_{u^2u^2}
\end{eqnarray*}

Now we may similarly proceed the study of second order equations.

The second order Lie equations arise from the conditions

\begin{eqnarray*}
\eta^1_{11} \vert ({\cal S}_2) = 0, \eta^1_{11} \vert ({\cal S}_2) = 0,
\eta^1_{11} \vert ({\cal S}_2) = 0,
\end{eqnarray*}

After direct computations we obtain the following groups of equations:

\begin{eqnarray*}
& &\eta^1_{x_1x_1} = 0,
\eta^1_{x_1x_2} = 0,
\eta^1_{x_2 x_2} = 0,\\
& &\theta^1_{u^1u^1} +
2a_{11}\theta^1_{u^1u^2} +
a^2_{11}\theta^1_{u^2u^2} = 0,-a_{21}\left (  \theta^1_{u^1 u^2} +
a_{11}\theta^1_{u^2u^2} \right ) = 0, a^2_{21}\theta^1_{u^2u^2} = 0,\\
& &\theta^2_{u^1u^1} +
2a_{22}\theta^2_{u^1 u^2} +
a^2_{22}\theta^2_{u^2u^2} = 0, -a_{12}\left ( \theta^2_{ u^1 u^2} +
a_{22}\theta^2_{u^2u^2} \right ) = 0, a^2_{12}\theta^2_{u^2u^2} = 0
\end{eqnarray*}

We have the following equations for $\eta^1$ and $\theta_1$:

\begin{eqnarray*}
& & 2 \eta^1_{ x_1  u^1} - \theta^1_{x_1x_1} +
2a_{11}\eta^1_{ x_1 u^2} = 0, 
\eta^1_{x_2 u^1} -
\theta^1_{x_1  x_2} +
a_{11}\eta^1_{x_2  u^2} +
a_{21}\eta^1_{x_1u^2} = 0,\\
& &\eta^1_{ u^1u^1} -
2\theta^1_{x_1  u^1}+
2a_{11}\left ( \eta^1_{ u^1 u^2} -
\theta^1_{x_1  u^2} \right ) +
a^2_{11}\eta^1_{ u^1u^1} = 0,\\
& &- \theta^1_{ x_2 u^1} -
a_{11}\theta^1_{x_2 u^2} +
a_{21}\left ( \eta^1_{u^1 u^2} -
\theta^1_{ x_1  u^2} \right ) +
a_{11}a_{21}\eta^1_{u^2u^2} = 0,\\
& &- \theta^1_{x_2x_2} +
2a_{21}\eta^1_{ x_2  u^2} = 0,
a_{21}\left ( - 2\theta^1_{ x_2 u^2} +
a_{21}\eta^1_{u^2u^2} \right ) = 0
\end{eqnarray*}

We also have similar equations for $\eta^1$, $\theta_2$:

\begin{eqnarray*}
& & 2 \eta^1_{ x_2  u^1} -
\theta^2_{x_2x_2} +
2a_{22}\eta^1_{ x_2 u^2} = 0,
\eta^1_{x_1 u^1} -
\theta^2_{x_1 x_2} +
a_{22}\eta^1_{ x_1 u^2} +
a_{12}\eta^1_{x_2 u^2} = 0,\\
& &\eta^1_{u^1u^1} -
2\theta^2_{x_2 u^1}+
2a_{22}\left ( \eta^1_{u^1 u^2} -
\theta^2_{ x_2  u^2} \right ) +
a^2_{22}\eta^1_{ u^2u^2} = 0,\\
& &- \theta^2_{x_1 u^1} -
a_{22}\theta^2_{x_1  u^2} +
a_{12}\left ( \eta^1_{u^1 u^2} -
\theta^2_{x_2  u^2} \right ) +
a_{22}a_{12}\eta^1_{u^2u^2} = 0,\\
& &- \theta^2_{x_1x_1} +
2a_{12}\eta^1_{x_1u^2} = 0,
a_{12}\left ( - 2\theta^2_{x_1 u^2} +
a_{12}\eta^1_{u^2u^2} \right ) = 0
\end{eqnarray*}

We have also the "mixed" equations containing $\eta^1$ and
both of $\theta_1$, $\theta_2$:

\begin{eqnarray*}
& &-2 \theta^2_{x_1 u^1}
+ 2 a_{12} \left ( \eta^1_{ u^1 u^2}
- \theta^1_{x_1 u^2}  \right ) -
2a_{11}\theta^2_{ x_1  u^2} +
2a_{11}a_{12}\eta^1_{ u^2u^2} = 0,\\
& & \eta^1_{u^1u^1} -
\theta^1_{ x_1  u^1}  -
\theta^2_{x_2  u^1}  -
a_{12}\theta^1_{ x_2 u^2} +
a_{11}\left ( \eta^1_{u^1  u^2} -
\theta^2_{x_2 u^2} \right )\\
& &+ a_{22} \left ( \eta^1_{ u^1 u^2} -
\theta^1_{x_1  u^2} \right ) -
a_{21}\theta^2_{ x_1  u^2} +
(a_{12}a_{21} + a_{11}a_{22})\eta^1_{u^2u^2} = 0,\\
& &-\theta^1_{x_2 u^1} -
a_{22}\theta^1_{x_2  u^2} +
a_{21} \left ( \eta^1_{u^1 u^2} -
\theta^2_{x_2  u^2} \right ) +
a_{21}a_{22}\eta^1_{u^2u^2} = 0
\end{eqnarray*}

Finally, we have the following series of equations :

\begin{eqnarray*}
& &\theta^2_{u^1 u^1} +
2a_{12}\theta^1_{ u^1  u^2} +
2a_{11}\theta^2_{ u^1  u^2} +
2a_{11}a_{22}\theta^1_{u^2u^2} +
a^2_{11}\theta^2_{u^2u^2} = 0,\\
& &a_{12}\left ( 2\theta^2_{u^1 u^2}
+ 2a_{11}\theta^2_{u^2u^2} +
a_{12}\theta^1_{u^2u^2} \right ) = 0,\\
& &\theta^1_{u^1u^1}  + (a_{11} + a_{22})
\theta^1_{ u^1 u^2} +
a_{21}\theta^2_{u^1  u^2} +
a_{11}a_{21}\theta^2_{ u^2u^2} +
(a_{11}a_{22} + a_{12}a_{21})\theta^1_{u^2u^2} = 0,\\
& &\theta^2_{ u^1u^1} + (a_{11} + a_{22})
\theta^2_{u^1 u^2} +
a_{12}\theta^1_{u^1  u^2} +
a_{12}a_{22}\theta^1_{u^2u^2}
+ (a_{11}a_{22} + a_{12}a_{21})\theta^2_{u^2u^2} = 0,\\
& &a_{21}\left ( 2\theta^1_{u^1  u^2} +
a_{21}\theta^2_{ u^2u^2} +
2a_{22}\theta^1_{u^2u^2} \right ) = 0,\\
& &\theta^1_{ u^1u^1} +
2a_{22}\theta^1_{ u^1 u^2}
+ 2a_{21}\theta^2_{u^1  u^2} +
2a_{21}a_{22}\theta^2_{u^2u^2} +
a^2_{22}\theta^1_{u^2u^2} = 0
\end{eqnarray*}

These equations together with earlier obtained first order ones form
the system of Lie equations for $({\cal S})$.

In order to  show  that the obtained second order linear PDE system is
of finite type and the type is equal to 1 it is necessary  to  study 
the 1-prolongation of this system i.e. essentially  the PDE system
obtained
by the consideration the first order partial derivatives of our equations.

  Two cases can occur: the case where  $a_{12} \neq 0$ or $a_{21} \neq 0$ 
and the case where  $a_{12} = a_{21} = 0$ and $a_{11} \neq a_{22}$. 
In every case the direct elementary computation shows that
 the symbol of the 1-prolongation is trivial.

This completes the proof of the proposition.

As a corollary we obtain the following

\begin{e-cor}
Let  $({\cal S}^{\varepsilon})$:

\begin{eqnarray*}
& & u^j_{x_1x_1} = F^j_{11}(\varepsilon,x,u,u_x),
u^j_{x_1x_2} = F^j_{12}(\varepsilon,x,u,u_x),
u^j_{x_2x_2} = F^j_{22}(\varepsilon,x,u,u_x), j = 1,2\\
& & u^2_{x_1} = a_{11}u^{1}_{x_1} + a_{12}u^{1}_{x_2}
+G_1(\varepsilon,x,u,u^1_x)\\
& & u^2_{x_2} = a_{21}u^{1}_{x_1} + a_{22}u^1_{x_2} +
G_2(\varepsilon,x,u,u^1_x)\\
\end{eqnarray*}
be a holomorphic completely integrable deformation  of the
flat nondegenerate  system $({\cal S}^0) = ({\cal S})$. Then for every $\varepsilon$ close
to
the origin
enough one has $dim Lie({\cal S}^{\varepsilon}) \leq dim Lie({\cal S}^0)$
and
every inifinitesimal symmetry of $({\cal S}^{\varepsilon})$ is determined by
its second order Taylor expansion at the origin.
\end{e-cor}

In particular, since the Segre family of a 6-dimensional real 
analytic Levi-nodegenerate manifold in $\cc^4$ is decribed by a system of this
class, the present method allows to obtain explicit recurcive formulae 
for infinitesimal automorphisms of such a manifold.

In conclusion of this paper we emphasize again that our method can be used in 
order to obtain a very precise information on  automorphisms of wide 
classes of CR manifolds and related PDE systems. For instance, if we
replace the condition (i) in the definition of a Levi nondegenerate
manifold by the slightly weaker condition of the triviality of the kernel of
the Levi form, the Segre family will be given by a ``mixed'' PDE
system containing second order partial derivatives of several
dependent 
variables and first order equations with linear parts satisfying some 
independence conditions; our method works for this class of systems
with minor modifications. The condition (ii) of the Levi nondegeneracy 
also can be replaced by a weaker assumptions on the highest Levi
forms. This leads to systems where the terms of highest order (in the 
first order equations) satisfy some independence conditions.    The most powerful algebraic
tool for the study of the related Lie equations is given by the Spencer cohomology theory 
and the Cartan - Kahler theory of normal forms of analytic linear PDE systems 
(see \cite{Po}). Finally, the consideration of manifolds with the
degenerate first Levi form leads to PDE systems which are not solved with respect to the highest partial derivatives. The study of their Lie symmetries 
needs more advanced tools of the local complex analytic geometry.
Our approach also raises several other natural questions: equivalence 
problems and invariants of involutive second order PDE systems with 
first order relations, classifications of these systems with respect to 
the properties of symmetry group (non-compact, transitive, etc.) by 
analogy with very well known result of geometric complex analysis. But
perhaps the most important problem is to develop in a systematic way the
geometry of the Segre families of  real analytic CR manifolds from the
complex differential and algebraic geometry standpoint.

{\small Univesrit\'e des Sciences et Technologies de Lille, Laboratoire
d'Arithm\'etique
 - G\'eom\'etrie - Analyse - Topologie, Unit\'e Mixte de Recherche 8524,
U.F.R. de
Math\'ematique, 59655 Villeneuve d'Ascq, Cedex, France}


\begin{thebibliography}{99}

{\footnotesize


\bibitem{BER1} M.S.Baouendi, P.Ebenfelt, L.P.Rothschild, {\it 
Algebraicity of holomorphic mappings between real algebraic sets in 
$\cc^n$}, Acta Math. {\bf 177}(1996), 225-273.

\bibitem{BER2} M.S.Baouendi, P.Ebenfelt, L.P.Rothschild, {\it Rational 
dependence of smooth and analytic CR mappings on their jets},
Math. Ann. {\bf 315} (1999), 205-249. 

\bibitem{BER3} M.S.Baouendi, P.Ebenfelt, L.P.Rothschild, {\it CR automorphisms
of real analytic CR manifolds in complex space},
Comm. Anal. Geom. {\bf 6} (1998), 291-315.




\bibitem{BP2} E.Bedford, S.Pinchuk, {\it Convex domains with noncompact
automorphism groups}, Mat. Sb. {\bf 185} (1994), 3-26.

\bibitem{Be1} V.Beloshapka, {\it A uniqueness theorem for
automorphisms of a nondegenerate surface in a complex space},
Math. Notes
{\bf 47} (1990), 239-242.


\bibitem{Be2} V.Beloshapka, {\it On holomorphic transformations of
quadric}, Math. USSR Sb. {\bf 72} (1992), 189-205.




\bibitem{BlKu} G.W.Bluman, S.Kumei,  {\it Symmetries and differential
equations},
Springer-Verlag, 1989.

\bibitem{BM} S.Bohner, J.Martin, {\it Several complex variables}, 
Princeton Univ. Press, 1948.

\bibitem{Ca} E.Cartan, {\it Sur la g\'eom\'etrie pseudoconforme des
hypersurfaces de
deux variables complexes}, Ann. Math. Pura Apll. {\bf 11} (1932)17-90.

\bibitem{CM} S.S.Chern, J.K.Moser, {\it Real hypersurfaces in complex
manifolds},
Acta Math. {\bf 133} (1974), 219-271.

\bibitem{Ch} S.S.Chern, {\it On the projective structure of a real
hypersurface in
$\cc^{n+1}$}, Math. Scand. {\bf 36} (1975), 74-82.

\bibitem{Di} L.E.Dickson, {\it Differential equations from the group
standpoint,} Ann. Math. {\bf 25}(1924), 287-378.

\bibitem{DW} K.Diederich, S.Webster, {\it A reflection principle for
degenerate
real hypersurfaces}, Duke Math. J. {\bf 47} (1980), 835-845.

\bibitem{DF} K.Diederich, J.E.Fornaess, {\it Proper holomorphic mappings
between
real analytic pseudoconvex domains in $\cc^n$}, Math. Ann. {\bf 282} (1988),
681- 700.

\bibitem{DP} K.Diederich, S.Pinchuk, {\it Proper holomorphic maps in
dimension 2 extend},
Indiana Univ. Math. J. {\bf 44} (1995), 1089-1126.



\bibitem{EIS} V.Ezhov, A.Isaev. G.Schmalz, {\it Invariants of elliptic
and hyperbolic CR-structures of codimension 2},
Internat. J. Math. {\bf 10} (1999), 1-52.


\bibitem{Fa} J.Faran, {\it Segre families and real hypersurfaces}, 
{\bf 60}(1980), 135-172.


\bibitem{GL} F. Gonzalez-Gascon, A.Gonzalez-Lopez, {\it Symmetries of 
differential equations. IV}, J. Math. Phys. {\bf 24}(1983), 2006-2021.


\bibitem{LS} S.Lie, G.Scheffers, {\it Vorlesungen uber Continuierliche 
Gruppen}, Chelsea, Bronx, NY. 1971.


\bibitem{Lo} A.Loboda, {\it Real analytic generating manifolds of 
codimension 2 in $\cc^4$ and their biholomorphic mappings}, Math. USSR 
Izv. {\bf 33} (1989), 295-315.


\bibitem{Ol} P.Olver, {\it Applications of Lie Groups to differential
equations},
Springer-Verlag, 1986.

\bibitem{Ol2} P.Olver, {\it Equivalence, invariants and symmetry}, 
Cambridge Univ. Press. 1995.
\bibitem{Ov} L.V.Ovsiannikov, {\it Group Analysis of Differential
equations},
Academic Press, New York, 1982.


\bibitem{Po} J.-F.Pommaret, {\it Systems of partial differential equations
and Lie pseudogroups}, Gordon and Breach Sci. Publ. 1978.

\bibitem{Po2} J.-F.Pommaret, {\it Partial differential equations and 
group theory}, Kluwer, 1994.

\bibitem{Se} B.Segre, {\it Intorno al problem di Poincar\'e della
representazione pseudo-conform}, Rend. Acc. Lincei, {\bf 13} (1931),
676-683.

\bibitem{Sh} B.Shiffman, {\it Projective geometry and Poincare's
theorem 
on automorphisms of the ball}, Enseign. Math. {\bf 41} (1995), 201-216.

\bibitem{St} N.Stanton, {\it Infinitesimal CR automorphisms of real 
hypersurfaces}, Amer. J. Math. {\bf 118}(1996), 209-233.

\bibitem{Su} A.Sukhov, {\it Segre varieties and Lie symmetries}, Pub. IRMA,
Lille, 1999, V.50.

\bibitem{Su2} A.Sukhov, {\it On CR mappings of real quadric
manifolds},
Mich. Math. J. {\bf 41}(1999), 143-150.

\bibitem{Tr} A.Tresse, {\it D\'etermination des invariants ponctuels de
l'\'equation
differentielle du second ordre $y'' = \omega(x,y,y')$}, Hirzel, Leiptzig,
1896.

\bibitem{Ta} N.Tanaka, {\it On the pseudo-conformal geometry of
hypersurfaces of the
space of $n$ complex variables}, J. Math. Soc. Japan, {\bf 14} (1962),
397-429.


\bibitem{Tu} A.Tumanov, {\it Finite-dimensionality of the group of CR 
automorphisms of a standard CR manifold and proper holomorphic
mappings of Siegel domains} Funkts. Anal. i Pril. {\bf 17}(1983), 49-61.


\bibitem{We1} S.M.Webster, {\it On the Mapping Problem for algebraic real
hypersurfaces} Inv. Math. {\bf 43} (1977), 53-68.

\bibitem{We2} S.M.Webster, {\it On the reflection principle in several 
complex variables}, Proc. AMS {\bf 71}(1978), 26-28.

\bibitem{We3} S.M.Webster, {\it Double valued reflection in the
complex 
plane}, Enseign. Math. {\bf 42}(1996), 25-48.

\bibitem{We4} S.Webster, {\it Some birational invariants for algebraic
real
hypersurfaces}, Duke Math. J. {\bf 45}(1978), 39-46.

\bibitem{Za} D.Zaitsev, {\it Germs of local automorphisms of
real-analytic CR structures and analytic dependence of k-jets}, 
Math. Res. Letters {\bf 41}(1997), 823-842.

}
\end{thebibliography}
\end{document}